\documentclass[11pt]{article}
\usepackage{color}
\usepackage{amsmath}
\usepackage[width=7.7in, height=9in]{geometry}
\catcode`\@=11 \@addtoreset{equation}{section}
\def\theequation{\thesection.\arabic{equation}}
\catcode`\@=12

\usepackage{amssymb}
\usepackage{amsfonts}
\usepackage{color}
\usepackage{amssymb}
\usepackage{amsbsy} 

\DeclareMathOperator*{\esssup}{ess\,sup}

\title{\textbf{The incompressible limit of compressible finitely extensible nonlinear bead-spring chain models for dilute polymeric fluids}}
\author{\sc{Endre S\"uli and Aneta Wr\'oblewska-Kami\'nska}}

\date{\today}

\begin{document}
\maketitle \maketitle
\newtheorem{definition}{Definition}[section]
\newtheorem{theorem}{Theorem}[section]
\newtheorem{lemma}[theorem]{Lemma}
\newtheorem{corollary}[theorem]{Corollary}
\newtheorem{proposition}[theorem]{Proposition}
\newtheorem{remark}{Remark}[section]
\renewcommand{\theequation}{\thesection.\arabic{equation}}
\newcommand{\nc}{\newcommand}
\nc{\R}{{\mathbb R}}
\nc{\N}{{\mathbb N}}
\nc{\qed}{\mbox{}\nolinebreak\hfill \rule{2mm}{2mm}}
\nc{\weak}{\rightharpoonup}
\nc{\weakstar}{\stackrel{\ast}{\rightharpoonup}}
\nc{\proof}{{\bf Proof: }}

\newcommand{\vrho}{\varrho}
\nc{\modular}[1]{{\stackrel{ #1}{\longrightarrow\,}}}

\def\bbbone{{\mathchoice {\rm 1\mskip-4mu l}
{\rm 1\mskip-4mu l} {\rm 1\mskip-4.5mu l} {\rm 1\mskip-5mu l}}}

\def\tens#1{\pmb{\mathsf{#1}}}
\def\vec#1{\boldsymbol{#1}}

\newcommand{\vu}{\vec{u}}
\newcommand{\divx}      {{{\mathrm{div}}_x}\,}
\newcommand{\divq}      {{{\mathrm{div}}_q}\,}
\newcommand{\vr}		{\vrho}
\newcommand{\vre}		{\vr_\varepsilon}
\newcommand{\vrez}	    {\vr_{0,\varepsilon}}
\newcommand{\vrs} 		{\overline{\vr}}
\newcommand{\vret}		{\tilde{\vr}_\ep}
\newcommand{\re}		{\rho_\ep}
\newcommand{\rs}		{\bar\rho}
\newcommand{\ue}		{\vec{u}_\varepsilon}
\newcommand{\uez}		{\vec{u}_{0,\varepsilon}}
\newcommand{\ep}		{\varepsilon}
\newcommand{\tS}		{\tens{S}}
\newcommand{\q}		    {\vec{q}}
\newcommand{\n}		    {\vec{n}}
\newcommand{\vU}		{\vec{U}}
\newcommand{\ess} 	    {{\rm{ess}}}
\newcommand{\res}		{{\rm{res}}}
\newcommand{\dx}		{\,{\rm d}\vec{x}}
\newcommand{\dt}		{\, {\rm d}t}
\newcommand{\dq}		{\,{\rm d}\vec{q}}
\newcommand{\dxdt}	    {\, \dx\,{\rm d}t}
\newcommand{\len}		{\Delta_{\ep, {\rm N}}}
\newcommand{\He}		{\vec{H}_\ep}
\newcommand{\U}		    {\vec{U}}
\newcommand{\vphi}	    {\vec{\phi}}
\newcommand{\tD}		{\tens{D}}

\def\bbbone{{\mathchoice {\rm 1\mskip-4mu l}
{\rm 1\mskip-4mu l} {\rm 1\mskip-4.5mu l} {\rm 1\mskip-5mu l}}}

\date{} 

\abstract{
We explore the behaviour of global-in-time weak solutions to a class of bead-spring chain models, with finitely extensible nonlinear elastic (FENE) spring potentials, for dilute polymeric fluids. In the models under consideration the solvent is assumed to be a compressible, isentropic, viscous, isothermal Newtonian fluid, confined to a bounded open domain in $\R^3$, and the velocity field is assumed to satisfy a complete slip boundary condition. We show that as the Mach number tends to zero the system is driven to its incompressible counterpart.
}

\smallskip


\section{Introduction}
The purpose of this paper is to explore the singular limit of finite-energy global-in-time weak solutions, as the Mach number tends to zero,  to a class of kinetic models of dilute polymeric fluids with noninteracting polymer chains, where the solvent is a compressible, isentropic, viscous isothermal Newtonian fluid, and the polymer molecules suspended in the solvent are idealised as linear bead-spring chains with finitely extensible nonlinear elastic (FENE) spring potentials. The existence of finite-energy global weak solutions to this class of models, which are coupled compressible Navier--Stokes--Fokker--Planck systems, has been shown in \cite{BS2014} in the case of a nonslip boundary condition for the velocity of the solvent. With minor modifications, the existence proof presented in \cite{BS2014} extends to the case of a complete slip (also referred to as Navier slip) boundary condition and, for the sake of simplicity of the exposition, it is the latter complete slip boundary condition that we shall assume hereafter. By performing a rigorous passage to the limit we show that, as the Mach number tends to zero, the compressible Navier--Stokes--Fokker--Planck system is driven to its incompressible limit, which is the coupled FENE-type bead-spring chain model arising from the kinetic theory of dilute solutions of polymeric liquids with noninteracting polymer chains. The limiting model involves the unsteady incompressible Navier--Stokes equations, with an elastic extra-stress tensor appearing on the right-hand side of the momentum equation. The elastic extra-stress tensor stems from the random movement of the polymer chains and is defined by the Kramers expression through the associated probability density function that satisfies a Fokker--Planck-type parabolic equation. The existence of global-in-time weak solutions to the limiting incompressible model was shown in \cite{BS10} in the case of a nonslip boundary condition for the velocity of the solvent. Again, the existence proof presented in \cite{BS10} extends to the case of a complete slip boundary condition. Our main result is therefore that the incompressible Navier--Stokes--Fokker--Planck system considered in \cite{BS10}, but with a complete slip boundary condition, is the singular limit, as the Mach number tends to zero, of the compressible model considered in \cite{BS2014}, with a complete slip boundary condition.

Our proofs rely on combining the theoretical machinery developed in the monograph of Feireisl and Novotn\'y \cite{FN}, devoted to the rigorous analysis of various asymptotic limits of the compressible Navier--Stokes--Fourier system, with the analytical techniques involved in the existence proofs in \cite{BS10}, \cite{BS2014}, and \cite{FLS2016} for coupled incompressible and compressible Navier--Stokes--Fokker--Planck systems.

For fully macroscopic models of compressible non-Newtonian fluids there are relatively few results in the literature concerning the incompressible limit.
In \cite{LZ2005}, Lei and Zhou studied the existence of global classical solutions of the incompressible Oldroyd-B system on the two-dimensional torus. The Oldroyd-B model arises in the description of viscoelastic fluid flow and consists of a coupling of the Navier--Stokes equations with a system of first-order partial differential equations for the elastic extra stress tensor. Using well-posedness and stability results for the compressible case they pass to the incompressible limit. This strategy is based on dispersive energy estimates for the compressible Oldroyd-B model. The outcome of their analysis is a global existence/uniqueness result for solutions to the two-dimensional incompressible Oldroyd-B model for small data and a justification of the incompressible Oldroyd-B model as a limit of the slightly compressible Oldroyd-B model. In \cite{L2006},
Lei studied the incompressible limit problem for the compressible Oldroyd-B model on a torus, and showed that compressible flows with well-prepared initial data converge to incompressible flows as the Mach number tends to zero. The case of bounded domains $\Omega \subset \mathbb{R}^3$ was explored by Guillop\'{e}, Salloum and Talhouk \cite{GST2010}, where they proved the local and global existence of strong solutions to the weakly compressible Oldroyd-B model, and, with the aid of conformal coordinates, they also studied the incompressible limit problem with well-prepared initial data. Subsequently Fang and Zi \cite{FS2014} proved, using the techniques of Danchin \cite{D2007} in scale-invariant Besov spaces on $\mathbb{R}^d$, that solutions to the compressible Oldroyd-B model with ill-prepared initial data (in which case strong time-oscillations of solutions need to be considered) converge to those of the corresponding incompressible Oldroyd-B model as the Mach number tends to zero. More recently, Ren and Ou \cite{RO2015} showed the local existence of strong solutions to an Oldroyd-B model for incompressible viscoelastic fluids in a bounded domain
$\Omega \subset \mathbb{R}^d$, $d=2,3$,  via the incompressible limit: the main idea of the paper was to derive uniform estimates with respect to the Mach number for a linearised compressible Oldroyd-B system.

To the best of our knowledge no rigorous asymptotic results are available in the literature for the incompressible limit of global weak solutions (with large initial data) to kinetic models of dilute compressible polymeric fluids involving the (compressible) Navier--Stokes--Fokker--Planck system, and it is the study of this question that is our objective here. The paper is structured as follows: in the next two subsections of this introductory section we perform a nondimensionalisation of the compressible Navier--Stokes--Fokker--Planck system and we state our main result; Section \ref{sec2} is then devoted to the proof of our main result. We end the paper, in Section \ref{sec3}, with concluding remarks concerning possible extensions.

\subsection{Dimensionless form of the Navier--Stokes--Fokker--Planck system}

We consider the following system of nonlinear partial differential equations consisting of the equations of continuity and balance of linear momentum, having the form of the compressible Navier--Stokes equations in which the elastic extra-stress tensor (the polymeric part of the Cauchy stress tensor) appears as a source term in the balance of linear momentum equation. In particular, for a given $T \in \R_{>0}$ and a bounded open domain $\Omega \subset \R^3$ with $\partial\Omega \in C^{2,\alpha}$, $\alpha \in (0,1)$, we denote by $\rho: \Omega \times [0,T] \to \R$ the density of the solvent, by $\vu: \overline{\Omega} \times [0,T] \to \R^3$ the velocity of the solvent, which satisfies the conservation of linear momentum equation for a viscous compressible fluid, and $\psi:\overline\Omega\times \overline D \times [0,T] \to \R_{\geq 0}$ is the probability density function associated with the random movement of the polymer chains and satisfying a Fokker--Planck equation to be stated below. The configuration space domain $D:=D_1 \times \cdots \times D_K \subset \mathbb{R}^{3K}$ is the Cartesian product of $K$ balanced convex open sets $D_i \subset \mathbb{R}^d$, $i=1,\dots,K$, with $q_i \in D_i$ denoting the orientation vector of the $i$-th spring in the bead-spring chain representing an idealisation of a polymer chain suspended in the solvent. The coupled Navier--Stokes--Fokker--Planck system under consideration has the following form:
	\begin{equation}\label{NSFP_1}
	\begin{split}
	\partial_t \rho +  \divx (\rho \vu)  & =  0 \quad \mbox{ in } \Omega \times (0,T]
	, \\
	\partial_t (\rho \vu ) + \divx (\rho \vu \otimes \vu ) + \frac{1}{{\rm Ma}^2} \nabla_x p(\rho)
	& =     \frac{1}{\rm Re}\, \divx \tS (\vu)  + \frac{1}{{\rm Fr}^2}\,  \rho \vec{f}  + \frac{1}{\rm Re}\,\divx \tens{\tau}_1  - \frac{\tilde\xi}{{ \rm Ma}^2} \nabla_x \varrho^2
	\quad \mbox{ in } \Omega \times (0,T]
	, \\
	\partial_t \psi + \divx (\vu \psi)
	+ \sum_{i=1}^K {\rm div}_{\q_i} ((\nabla_x u) \q_i \psi) & =  \delta \Delta_x \psi
	+ \frac{1}{4 {\mathrm{De}}} \sum_{i=1}^{K} \sum_{j=1}^{K} A_{ij}\, {\rm div}_{\q_i}
	\left( M \nabla_{\q_j} \left(\frac{\psi}{M}\right)\right)  \quad \mbox{ in } \Omega \times D \times (0,T],\\
	\end{split}
	\end{equation}
with the initial conditions
	\begin{equation}\label{in_data_1}
	\begin{split}
	\rho (\cdot,0)  & =  \rho_0(\cdot) \geq 0 \quad \mbox{ in } \Omega, \\
	(\rho \vu )(\cdot,0)  & =  (\rho_0 \vu_0 )(\cdot)   \quad \mbox{ in } \Omega, \\
	\psi(\cdot,\cdot,0) & =  \psi_0(\cdot,\cdot) \geq 0\quad \mbox{ on } \Omega \times D.
	\end{split}
	\end{equation}
In the above system  $\mathrm{Ma}$ is the Mach number, $\mathrm{Re}$ is the Reynolds number, $\mathrm{Fr}$ is the Froude number, $\mathrm{De}$ is the Deborah number, $\tilde{\xi}$ and $\delta$ are positive real numbers, $A_{ij}$, $i,j=1,\dots,K$, in the Fokker--Planck equation \eqref{NSFP_1}$_3$, are the entries of 
\begin{equation}\label{as_AA}
\tens{A} = [A_{ij}]_{i,j=1}^K,  \quad \mbox{a symmetric positive definite matrix, with smallest eigenvalue }a_0\in \R_{>0},
\end{equation}
called the Rouse matrix, or connectivity matrix (e.g. $\tens{A} = {\rm tridiag}[-1,2,-1]$ in the case of a topologically linear bead-spring chain). The function $p=p(\rho)$ is the pressure, $\tS$ is the Newtonian part of the viscous stress tensor,  defined by
	\begin{equation}\label{tensS}
	\tS(\vu) := \mu^S \bigg( \tD \vu - \frac{1}{3} (\divx \vu)\, \tens{I}\bigg) + \mu^B (\divx \vu)\, \tens{I} ,
	\end{equation}
where $\tens{I}$ is the $3 \times 3$ identity tensor,
$\tD \vec{v}:= \frac{1}{2} (\nabla_x \vec{v} + (\nabla_x \vec{v})^{\rm T})$ is the rate of strain tensor, $\mu^S$ and $\mu^B$ are positive constants referred to as the shear viscosity and the bulk viscosity, respectively, and  the function $M=M(\vec{q})$ is the Maxwellian, to be definedhttps://www.overleaf.com/project/5bd0b4be1ae53d317f17aefd below, associated with the model.
In addition, the pressure $p$ is assumed to be related to the density $\rho$ of the solvent by the isentropic equation of state:
	\begin{equation}\label{p1}
	p(\rho) := c_p \rho^\gamma,
	\end{equation}
where $c_p \in \R_{>0}$ and the constant $\gamma > 3/2$.
The above system is supplemented with boundary conditions, which will be stated below.
	
In a bead-spring chain model consisting of $K+1$ beads, linearly coupled with $K$ elastic springs to represent a polymer chain, the extra-stress tensor $\tens{\tau}$  is given by a version of the Kramers expression depending on the probability density function
$\psi$ of the random conformation vector  $\vec{q} = (\vec{q}^{\rm T}_1, \dots, \vec{q}^{\rm T}_K)^{\rm T} \in D:= D_1 \times \cdots \times D_K \subset \R^{3K}$ of the chain, where $\vec{q}_i$  represents the 3-component conformation/orientation vector of the $i$-th spring in the chain.
Typically $D_i$ is the whole of $\R^3$ or a bounded open $3$-dimensional ball centred at the origin $\vec{0} \in \R^3$, for each $i = 1, \dots, K.$ If $K =1$, then the model is referred to as the \textit{dumbbell model}.
Here we shall concentrate on finitely extensible nonlinear elastic (FENE) bead-spring chain models, with
\begin{equation}\label{setD}
D := B(\mathbf{0}, b^{\frac{1}{2}}_1)  \times \cdots \times B(\mathbf{0}, b^{\frac{1}{2}}_K), \mbox{ where } b_i > 0,\ i=1, \dots, K,
\ K \geq 1,
\end{equation}
 and
$B(\mathbf{0}, b^{\frac{1}{2}}_i)$ is a bounded open ball in $\R^3$ of radius $b^{\frac{1}{2}}_i$, centred at $\mathbf{0} \in \R^3$.

On the right-hand side of \eqref{NSFP_1}$_2$, the 3-component vector function $\vec{f}$ is the nondimensional density of body forces and the elastic extra-stress tensor is of the form:
	\begin{equation}\label{def:tau}
	\tens{\tau} (\psi)(\vec{x},t)  :=  \frac{1}{\rm Re}\tens{\tau}_1(\psi) (\vec{x},t)
	- \frac{1}{{\rm Ma}^2} \bigg( \int_{D\times D}  \gamma (\vec{q},\vec{q}') \psi (\vec{x},\vec{q},t) \psi (\vec{x},\vec{q}',t)  {\rm \, d} \vec{q} {\rm \, d} \vec{q}' \bigg)\, \tens{I}
	=  \frac{1}{\rm Re} \tens{\tau}_1(\vec{x},t) - \frac{\tilde\xi}{{\rm Ma}^2} \varrho^2(\vec{x},t)\, \tens{I},
	\end{equation}
where $\gamma: D \times D \to \R_{\geq 0}$ is a smooth, time-independent, $x$-independent and $\psi$-independent interaction kernel, which we take here for the sake of simplicity to be
%
	$$\gamma (\vec{q}, \vec{q}') = \tilde\xi, \quad \mbox{ where } \tilde\xi \in \R_{\geq 0}, $$
and the polymer number density is defined by
	\begin{equation}\label{pol_num_den}
	\varrho (\vec{x},t) := \int_D \psi (\vec{x},\vec{q},t) {\rm\,d}\vec{q}, 	\quad (\vec{x},t) \in \Omega \times [0,T].
	\end{equation}
Moreover,  $\tens{\tau}_1(\psi)$ appearing in \eqref{def:tau} is  the Kramers expression:
	\begin{equation}\label{def_t1}
	\tens{\tau}_1(\psi) :=  \frac{1-\beta}{\mathrm{De}}
	\left(
	\sum_{i=1}^K \tens{C}_i ( \psi ) - (K+1) \int _D \psi {\rm\,d}\vec{q} \,\tens{I}
	\right),
	\end{equation}
where $(1-\beta) := \eta_p / (\eta_s + \eta_p)$,  with $\eta_p$ signifying the polymeric viscosity and $\eta_s$ the viscosity of the solvent. In the above
	\begin{equation}\label{C1}
	\tens{C}_i ( \psi ) (\vec{x},t) := \int_{D} \psi (\vec{x},\vec{q},t)\, U'_i\bigg(\frac{1}{2} | \vec{q}_i|^2\bigg) \vec{q}_i \vec{q}^{\rm T}_i {\rm\,d}\vec{q}, \quad
	i=1,\dots,K,
	\end{equation}
where $U_i$ is the $i$-th spring potential, which we shall now fix and we shall also define the associated Maxwellian.

Let ${\mathcal{O}}_i \subset [0,\infty)$ denote the image of $D_i$  ($0 \in {\mathcal{O}}_i$) under the mapping
$\vec{q}_i \in D_i \mapsto \frac{1}{2} |\vec{q}_i|^2$ and let us consider the spring potential $U_i \in C^1( {\mathcal{O}}_i; \R_{\geq 0} ),$ $i= 1,\dots,K$. In the case of the FENE spring potential, $\mathcal{O}_i = [0,b_i)$, $i=1,\dots, K$.
We shall suppose that $U_i(0) = 0$ and that $U_i$ is unbounded on ${\mathcal{O}}_i $ for $i= 1,\dots,K$ (i.e., $\lim_{s \rightarrow b_i{-}} U_i(s) = +\infty$). The elastic spring-force $\vec{F}_i : D_i \subseteq \R^3 \to \R^3$ of the $i$-th spring in the chain is then defined by
	\begin{equation}\label{defF}
	\vec{F}_i(\vec{q}_i) := U'_i\bigg( \frac{1}{2}|\vec{q}_i|^2 \bigg)\vec{q}_i, \quad i=1,\dots,K,
	\end{equation}
and the partial Maxwellian $M_i$, associated with the spring potential $U_i$, is defined by
	\begin{equation*}
	M_i(\vec{q}_i) := \frac{1}{Z_i} \mathrm{e}^{-U_i(\frac{1}{2} | \vec{q}_i |^2)}, \quad Z_i := \int_{D_i} {\mathrm e}^{-U_i(\frac{1}{2} | \vec{q}_i |^2)}
	{\rm d} \vec{q_i}, \quad i=1,\dots,K.
	\end{equation*}
The total Maxwellian $M_i$ in the model is then
	\begin{equation}\label{def_total_M}
	M(\vec{q}) := \prod_{i=1}^K M_i (\vec{q_i}) \quad \mbox{ for all } \vec{q}
	= (\vec{q}_1^{\rm T}, \dots, \vec{q}_K^{\rm T})^{\rm T} \in D = D_1 \times \cdots \times D_K.
	\end{equation}
A straightforward calculation yields that, for $i=1,\dots,K$,
	\begin{equation}\label{Max1}
	M(\vec{q}) \nabla_{q_i} [M(\vec{q})]^{-1} = - [ M(\vec{q})]^{-1} \nabla_{q_i} M(\vec{q})
	= U_i' \bigg(\frac{1}{2} |\vec{q}_i|^2\bigg) \vec{q}_i,
	\end{equation}
and
	\begin{equation}\label{Max2}
	\int_D M(\vec{q}) {\rm\,d}\vec{q} = 1.
	\end{equation}
We shall suppose that for $i=1,\dots,K$ there exist constants $c_{ij} >0,$ $j=1, \dots, 4$, and $\theta_i > 1$ such that the spring potential $U_i$ and the associated Maxwellian $M_i$ satisfy
	\begin{equation}\label{assU1}
	c_{i1} [{\rm dist} (\vec{q}_i, \partial D_i)]^{\theta_i} \leq M_i (\vec{q}_i)
	\leq c_{i2}  [{\rm dist} (\vec{q}_i, \partial D_i)]^{\theta_i} \quad \mbox{ for all } \vec{q}_i \in D_i,
	\end{equation}
	\begin{equation}\label{assU2}
	c_{i3} \leq  {\rm dist} (\vec{q}_i, \partial D_i)\, U'_{i}\bigg(\frac{1}{2} |\vec{q}_i|^2\bigg) \leq c_{i4} \quad \mbox{ for } \vec{q}_i \in D_i \mbox{ with } i=1,\dots, K.
	\end{equation}
It then follows from the above that
	\begin{equation}\label{U3}
	\int_{D_i} \left[ 1+ \bigg[U_i\bigg(\frac{1}{2} |\vec{q}_i|^2\bigg)\bigg]^2 + \bigg[U'_i\bigg(\frac{1}{2} |\vec{q}_i|^2\bigg)\bigg]^2    \right]
	M_i (\vec{q}_i) {\rm\,d}\vec{q}_i < \infty , \quad i=1, \dots, K.
	\end{equation}

In the system \eqref{NSFP_1} the Reynolds number Re and the Deborah number De are defined, respectively, by
	\begin{equation}
	{\rm Re} := \frac{\rho_0 L_0 U_0}{\eta_s + \eta_p},
	\end{equation}
	\begin{equation}
	{\rm De} := \frac{\lambda U_0}{L_0} = \frac{\zeta_0 U_0}{4 H L_0},
	\end{equation}
with $\mathrm{De}$ characterising the elastic relaxation properties, $\lambda>0$ being the characteristic relaxation time, $U_0$ the characteristic speed of the fluid, $L_0$ the characteristic macroscopic lengthscale (e.g. the diameter of the flow domain $\Omega$), $H>0$ the spring constant, $\zeta_0>0$ the characteristic drag coefficient, and $\rho_0$ is the characteristic density; the Mach number and the Froude number are defined, respectively, by
	 \begin{equation}
	{\rm Ma} := \frac{U_0}{ \sqrt{ p_0 / \rho_0} },
	\end{equation}
	\begin{equation}
	{\rm Fr} := \frac{U_0}{ \sqrt{ L_0 / f_0} },
	\end{equation}
where $p_0$ is the characteristic pressure, $f_0$ signifies the characteristic density of body forces; and
	$$\delta := \frac{(l_0 /  L_0)^2} {(4 (K+1) \zeta_0 U_0)/ (4H L_0)} = \frac{(l_0 /  L_0)^2} {4 (K+1) {\mathrm{De}}}$$
is the polymeric diffusion coefficient, where $l_0$ is the characteristic microscopic lengthscale (e.g. the length of a typical polymer chain suspended in the solvent).

We supplement the system of partial differential equations under consideration with a complete slip boundary condition for the velocity field of the solvent, that is,
\begin{equation}\label{bc-c-slip}
\vu \cdot \vec{n} |_{\partial \Omega} = 0, \quad\quad [\tS(\vu) + \tens{\tau}_1]\vec{n} \times \vec{n} |_{\partial \Omega} = \vec{0}.
\end{equation}
Moreover, we impose the following boundary and initial conditions on solutions of the Fokker--Planck equation:
	\begin{equation}\label{bc_psi_1}
	\left[
	\frac{1}{4 {\mathrm{De}}} \sum_{j=1}^K A_{ij} M \nabla_{q_j} \left( \frac{\psi}{M} \right) - (\nabla_x \vu) \vec{q}_i \psi
	\right]
	\cdot \frac{ \vec{q}_i }{| \vec{q}_i |} = 0 \quad \mbox{ on } \Omega \times \partial \overline{D}_i \times (0, T], \mbox{ for } i=1,\dots,K,
	\end{equation}
	\begin{equation}\label{bc_psi_2}
	\delta \nabla_x \psi \cdot \vec{n} = 0 \quad \mbox{ on } \partial \Omega \times D \times (0,T],
	\end{equation}
	\begin{equation}\label{bc_psi_3}
	\psi(\cdot,\cdot,0) = \psi_0 (\cdot, \cdot) \geq 0 \quad \mbox{ on } \Omega \times D,
	\end{equation}
where $\partial \overline{D}_i := D_1 \times \cdots \times D_{i-1} \times \partial D_i \times D_{i+1} \times \cdots \times D_K,$
$\vec{q}_i$ is normal to $\partial D_i$, as $D_i$ is a bounded ball centred at the origin, and $\vec{n}$ is the unit outward normal vector to $\partial \Omega.$

\begin{remark} The Navier--Stokes--Fokker--Planck system stated above is in dimensionless form. It is arrived at by introducing new variables of the form
$X^* = \frac{X}{X_0}$ into the dimensionful version of the system, where $X_0>0$ is the characteristic value of $X$; for example, $L_0>0$ is  a reference length, $T_0>0$ is a reference time, $U_0>0$ is a reference velocity, $\rho_0>0$ is a reference density of the solvent, proceeding similarly with the characteristic values of the other physical entities entering into the equations: $p_0>0$, $\mu^S_0>0$, $\mu^B_0>0$, $f_0>0$, $\eta^p>0$. In the above scaling the polymeric pressure term $p_p = \xi \vrho^2$ is treated as a part of a fluid pressure, and it is therefore rescaled by the reference value $p_0$.
\end{remark}

In this paper we concentrate on the following choices of the various dimensionless numbers:

\begin{itemize}
\item ${\rm Ma} = \ep \ll 1$, i.e., the characteristic velocity is dominated by the speed of sound;
\item ${\mathrm{De}} = 1$ and ${\rm Re} =1$;
\item One can choose ${\rm Fr} = \sqrt{\ep}$ (the case of low stratification), but for the sake of simplicity we shall not consider the influence of external forces in the momentum equation in this work, and will therefore set $\vec{f} = \vec{0}$;
\item For the dimensionless coefficient $\tilde\xi$ we shall assume here that $\tilde \xi = \bar\xi\, {\rm Ma}^2 = \bar\xi\, \ep^2$, with some $\bar\xi > 0$.
\end{itemize}

The system \eqref{NSFP_1} therefore takes the following form: for a given $T \in \R_{>0}$ and a bounded open domain $\Omega \subset \R^3$, with $\partial\Omega \in C^{2,\alpha}$, $\alpha \in (0,1)$, we consider the following system of partial differential equations:
	\begin{equation}\label{NSFP_2}\tag{NSFP$_\ep$}
	\begin{split}
	\partial_t \rho +  \divx (\rho \vu)  &  =   0 \quad \mbox{ in } \Omega \times (0,T]
	, \\
	\partial_t (\rho \vu ) + \divx (\rho \vu \otimes \vu ) + \frac{1}{\ep^2} \nabla_x p(\rho)
	&  =     \divx \tS (\vu)  + \divx \tens{\tau}_1  - {\bar\xi} \nabla_x \varrho^2
	\quad \mbox{ in } \Omega \times (0,T]
	, \\
	\partial_t \psi + \divx (\vu \psi)
	+ \sum_{i=1}^K {\rm div}_{\q_i} ((\nabla_x \vu) \q_i \psi) & =  \delta \Delta_x \psi
	+ \frac{1}{4} \sum_{i=1}^{K} \sum_{j=1}^{K} A_{ij} {\rm div}_{\q_i}
	\left( M \nabla_{\q_j} \bigg(\frac{\psi}{M}\bigg)\right) \quad \mbox{ in } \Omega\times D \times (0,T],
	\\ \mbox{ with } 	
	\varrho (\vec{x},t) & := \int_D \psi (\vec{x},\vec{q},t) {\rm\,d}\vec{q}, 	\quad (\vec{x},t) \in \Omega \times [0,T]
	\end{split}
	\end{equation}
for each $\ep \in (0,1),$
supplemented with the boundary conditions \eqref{bc-c-slip}--\eqref{bc_psi_2}, the initial conditions \eqref{in_data_1}, and the initial data appearing in \eqref{in_data_1} assumed to satisfy the following properties:
	\begin{equation}\label{indat_1}
	\| \vu_{0,\ep} \|_{L^2(\Omega)} \leq c,
	\end{equation}
	\begin{equation}\label{indat_2}
	\rho_{0,\ep} = \rs + \ep r_{\ep,0}
	\quad \mbox{ where } \quad
	\| r_{\ep,0} \|_{L^\infty(\Omega)} \leq c~
	 \quad  \mbox{ and } \quad \rs > 0,  \quad \rho_{0,\ep} \geq 0\quad \mbox{ a.e. } \vec{x}\in \Omega, \quad \int_{\Omega} r_{\ep,0} \,\mathrm{d}\vec{x} = 0,
	\end{equation}
	\begin{equation}\label{indat_3}
	\psi_{0,\ep} \geq 0 \mbox{ a.e. on }\Omega  \times D \quad  \mbox{ and } \quad \|{\mathcal{F}}(\widehat\psi_{0,\ep})\|_{L^1_{M}(\Omega \times D)} \leq c, \quad
	 \bigg\| \int_D \psi_{0,\ep} (\cdot,\vec{q}) {\rm\,d}\vec{q} \,\bigg\|_{L^{\infty}_{\geq 0} (\Omega)  } \leq c ,
	\end{equation}
where $\vrs$ is a static constant density satisfying the equilibrium state (static state) equation
	\begin{equation}\label{st_st}
	\nabla_x p(\rs) = 0 \quad \mbox{ in } \Omega.
	\end{equation}
In the above we have used the following notation:
$$\widehat\psi := \psi/M\qquad\mbox{and}\qquad
 \mathcal{F}(s)  := s (\log(s) - 1) \quad  \mbox{for $s > 0$},\quad \mbox{with $\mathcal{F}(0):=0.$}$$
Clearly,  $\mathcal{F}'(s) = \log(s)$ and $\mathcal{F}''(s) = \frac{1}{s}$ for $s>0$.

\begin{remark}[Remark~1.3 in \cite{BS2014}]
Defining the polymer number density by \eqref{pol_num_den}, formally integrating the Fokker--Planck equation in \eqref{NSFP_1} over $D$, and using the boundary condition \eqref{bc_psi_1}, we infer that
	\begin{equation}\label{rem_11}
	\partial_t \varrho + \divx (\varrho \vu) = \delta \Delta_x \varrho \quad \mbox{ on } \Omega \times (0,T] ,
	\end{equation}
with $\varrho(x,0) = \varrho_0$ in $\Omega$ and $\frac{\partial \varrho}{\partial n} |_{\partial \Omega} = 0$.
Hence by the boundary and initial conditions \eqref{bc_psi_2}, \eqref{bc_psi_3} we get that
	\begin{equation}\label{rem_12}
	\delta \nabla_x \varrho \cdot \vec{n} = 0 \quad \mbox{ on } \partial \Omega \times (0,T] \quad \mbox{ and } \quad
	\varrho (\vec{x},0) = \int_D \psi_0(\vec{x},\vec{q})\, {\rm d} \vec{q} \quad \mbox{ for } \vec{x} \in \Omega.
	\end{equation}
If $\divx  \vu = 0$ and $\varrho(\cdot, 0)$  is constant, then $\varrho(\vec{x},t)$ is constant  (and equal to $\varrho(\vec{x}, 0)\equiv \mathrm{Const.}>0$) for all $(\vec{x},t) \in \Omega\times(0,T]$, and thus 	
	$$ \varrho(\vec{x},t) =  \int_D \psi(\vec{x},\vec{q},t)\, {\rm d} \vec{q} = \int_D \psi_0(\vec{x},\vec{q})\, {\rm d} \vec{q}
	= \varrho(\vec{x},0)
	 \in \R_{>0} \quad  \mbox{ for all } (\vec{x},t) \in \Omega \times (0,T]. $$
In other words, the polymer number density is constant.
\end{remark}

Our aim is to show by rigorous analysis that, as the Mach number converges to zero, the primitive system \eqref{NSFP_2} converges to the following incompressible Navier--Stokes--Fokker--Planck system:
	\begin{equation}\label{limit1}
	\begin{split}
	\divx \vU  & =  0 \quad \mbox{ in } \Omega \times (0,T]
	, \\
	\partial_t (\rs \vU ) + \divx (\rs \vU \otimes \vU ) + \nabla_x \Pi
	 & =     \mu_S \divx \tD_x \vU  + \divx \tens{\tau}_1  (\Psi)
	\quad \mbox{ in } \Omega \times (0,T]
	, \\
	\partial_t \Psi + \divx (\vU \Psi)
	+ \sum_{i=1}^K {\rm div}_{\q_i} ((\nabla_x \vU) \q_i \Psi) & =  \delta \Delta_x \Psi
	+ \frac{1}{4} \sum_{i=1}^{K} \sum_{j=1}^{K} A_{ij} {\rm div}_{\q_i}
	\left( M \nabla_{\q_j} \left(\frac{\Psi}{M}\right)\right)  \quad \mbox{ in } \Omega \times (0,T],
	\end{split}
	\end{equation}
where $\Pi$ is a new pressure function,  $\vU$ is the limiting divergence-free velocity field, the density of the solvent $\rs$ is constant,
	\begin{equation}
	\tens{\tau}_1 :=  (1-\beta)
	\left(
	\sum_{i=1}^K \tens{C}_i ( \Psi ) - (K+1) \int _D \Psi {\rm\,d}\vec{q}\, \tens{I}
	\right)
	\end{equation}
is the Kramers expression for the polymeric extra stress tensor in the incompressible limit, with, again,  $(1-\beta) := \eta_p / (\eta_s + \eta_p)$, and
	\begin{equation}\label{C1_limit}
	\tens{C}_i ( \Psi ) (\vec{x},t) := \int_{D} \Psi (\vec{x},\vec{q},t)\, U'_i\bigg(\frac{1}{2} | \vec{q}_i|^2\bigg) \vec{q}_i \vec{q}^T_i {\rm\,d}\vec{q}, \quad
	i=1,\dots,K.
	\end{equation}

\begin{remark}
One can formally arrive at the system \eqref{limit1} by inserting in the primitive system the expansions:
	\begin{equation*}
	\begin{split}
	\rho & = \bar\rho + \ep \rho^{(1)} + \ep^2 \rho^{(2)} + \cdots, \\
	u &= \vec{U} + \ep \vu^{(1)} + \ep^2\vu^{(2)} + \cdots,\\
	\psi & = \Psi + \ep \psi^{(1)} + \ep^2 \psi^{(2)} + \cdots,
	\end{split}
	\end{equation*}
and neglecting all terms of order $\ep$ and higher. Our goal here is to deduce \eqref{limit1} from the primitive system \eqref{NSFP_2} through a rigorous passage to the limit, as $\ep \rightarrow 0$ (i.e, as $\mathrm{Ma} \rightarrow 0$).
\end{remark}

\subsection{Main result}
For $r \in[1,\infty)$, let  $L^r_M(\Omega \times D)$ denote the Maxwellian-weighted Lebesgue space equipped with the norm
	$$\| \varphi \|_{L^r_M} := \left( \int_{\Omega \times D} M |\varphi|^r {\rm\,d}\vec{q}\dx  \right)^{\frac{1}{r}},$$
and let $H^1_M(\Omega \times D)$ denote the Maxwellian-weighted Sobolev space with the norm
	$$\| \varphi \|_{H^1_M} := \left( \int_{\Omega \times D} M \left[  |\varphi|^2   + | \nabla_x \varphi |^2 + | \nabla_q \varphi |^2 \right]{\rm\,d}\vec{q}\dx  \right)^{\frac{1}{2}}.$$
We also define
	\begin{equation}\label{d:zm}
	Z_r := \{ \phi \in L^r_M(\Omega \times D)\,:\, \phi\geq 0\, \mbox{ a.e. on }\Omega \times D \}.
	\end{equation}

\begin{definition}\label{weak_sol_primitiv}
We call a  triple $(\rho_\ep, \vu_\ep, \widehat\psi_\ep)$ a weak solution to the system  \eqref{NSFP_2}
provided that:
\begin{itemize}
\item[1)] The functions $\rho_\ep$, $\vu_\ep$, $\widehat\psi_\ep$ satisfy the following regularity properties:
    \begin{equation}\label{L61_0}
	\rho_\ep \in C_{w} ([0,T]; L^\gamma_{\geq 0} (\Omega)) \cap H^1 (0,T; W^{1,6}(\Omega)') \cap  L^2(0,T; H^1(\Omega)'),
	\end{equation}
	\begin{equation}\label{L61_1}
	\vu_\ep \in L^2 (0,T ; H^1 (\Omega)^3) \quad \mbox{ and } \quad
	\widehat\psi_\ep \in L^{v}(0,T;Z_1) \cap H^1(0,T; M^{-1} (H^s(\Omega \times D))'),
	\end{equation}
 where $v \in [1,\infty)$ and $s > 1+ \frac{3}{2} (K +1) $, with finite relative entropy and Fisher information,
 	\begin{equation}\label{L61_2}
	{\mathcal{F}}(\widehat\psi_\ep) \in L^\infty(0,T; L^1_M(\Omega\times D)) \quad
	\mbox{ and } \sqrt{\widehat\psi_\ep} \in L^2(0,T; H^1_M(\Omega \times D)),
	\end{equation}
	\begin{equation}\label{L61_8}
	\tens{\tau}_1 (M \widehat\psi_\ep)  \in L^r(\Omega\times [0,T))^{3 \times 3}, \quad \mbox{ where } r \in \bigg[ 1, \frac{20}{13}  \bigg);
	\end{equation}
\item[2)] In addition,
	\begin{equation}\label{L61_10}
	\varrho_\ep := \int_D M\widehat\psi_\ep {\rm\,d}\vec{q} \in L^\infty(0,T;L^2(\Omega)) \cap L^2(0,T; H^1(\Omega)),
	\end{equation}
and
	\begin{equation}\label{L61_12}
	\varrho_\ep \in L^{\frac{5\zeta}{3(\zeta-1)}}(0,T; L^\zeta(\Omega )) \quad\quad \mbox{ for any }\zeta\in(1,6);
	\end{equation}
\item[3)] Moreover the following relations are satisfied:
 	\begin{equation}\label{T61_1}
	\int_0^T \left\langle  \partial_t \rho_\ep , \eta  \right\rangle_{W^{1,6}(\Omega)} \dt
	- \int_0^T \int_\Omega \rho_\ep \vu_\ep \cdot \nabla_x \eta \dx\dt = 0 	
	\quad \mbox{ for all } \eta \in L^2(0,T; W^{1,s}(\Omega)) ,
	\end{equation}
for some $s\in \left( \left(\frac{6\gamma}{6+\gamma}\right)',\infty \right]$, where $\left(\frac{6\gamma}{6+\gamma}\right)'$ is the H\"older conjugate of $\frac{6\gamma}{6+\gamma}$, and  with $\rho_\ep( \cdot, 0) = \rho_{0,\ep}$,
	\begin{equation}\label{T61_2}
	\begin{split}
	& \int_0^T \left\langle \partial_t(\rho_\ep \vu_\ep), \vec{w} \right\rangle_{W^{1,r}_0(\Omega)} \dt
	+ \int_0^T \int_\Omega \left[ \tens{S}(\vu_\ep) - \rho_\ep \vu_\ep \otimes \vu_\ep - \frac{1}{\ep}c_p \rho_\ep^\gamma \, \tens{I} \right] : \nabla_x \vec{w} \dx\dt
	\\ &
	= - \int_0^T \int_\Omega  \left( \tens{\tau}_1(M \widehat\psi_\ep) - \bar \xi \varrho_\ep^2\, \tens{I} \right)
	: \nabla_x \vec w
	\dx\dt
	  \quad \mbox{ for all } \vec{w} \in L^{\frac{\gamma + \vartheta}{\vartheta}}(0,T; W^{1,r}(\Omega)^3), \mbox{ where } \vec{w}\cdot \vec{n}|_{\partial\Omega}=0,
\end{split}
	\end{equation}
with $(\rho \vu)(\cdot,0) = (\rho_0\vu_0)(\cdot),$ $\vartheta(\gamma):= \frac{2\gamma - 3}{3}$ for $\frac{3}{2}<  \gamma \leq 4$
and $\vartheta(\gamma):= \frac{5}{12}\gamma$ for
$4\leq \gamma$; $r:= \max \left\{  4, \frac{6\gamma}{2\gamma - 3} \right\}$,
and
	\begin{equation}\label{T61_3}
	\begin{split}
	& \int_0^T \left\langle  M\partial_t \widehat\psi_\ep , \varphi   \right\rangle_{H^s (\Omega \times D)} \dt
	+  \frac{1}{4}   \sum_{i,j=1}^K A_{ij} \int_{\Omega\times D}  M \nabla_{q_j}  {\widehat\psi_\ep}
	\cdot  \nabla_{q_i}  {\varphi} {\rm\,d}\vec{q}\dx\dt
	\\ &
	\quad + \int_0^T \int_{\Omega \times D} M [ \delta \nabla_x \widehat\psi_\ep - \vu_\ep\widehat\psi_\ep ] \cdot \nabla_x \varphi {\rm\,d}\vec{q}\dx\dt
	- \int_0^T \int_{\Omega\times D }  M \sum_{i=1}^K [( \nabla_x \vu_\ep) \vec{q}_i ] \widehat\psi_\ep  \cdot \nabla_{q_i} \varphi
	{\rm\,d}\vec{q}\dx\dt
	\\ & = 0 \quad \mbox{ for all } \varphi \in L^2(0,T ; H^s(\Omega \times D)),
	\end{split}
	\end{equation}
with $\widehat\psi_\ep (\cdot, 0) = \widehat\psi_{0,\ep}(\cdot)$ and $s > 1 + \frac{3}{2}(K+1).$
\end{itemize}
\end{definition}

Let $\vec{H}$ denote the standard Helmholtz projection onto the space of solenoidal (divergence-free) functions; that is,
	$$
	\vec{v} = {\vec{H}}[\vec{v}] + {\vec{H}}^{\bot}[\vec{v}], \quad \vec{H}^{\bot}[\vec{v}] := \nabla_x \Phi,
	$$
where $\Phi \in H^1(\Omega)/\R$ is the unique solution of the problem
 	$$
	\int_\Omega \nabla_x \Phi \cdot \nabla_x \phi \dx = \int_\Omega \vec{v} \cdot \nabla_x \phi \dx
	\quad \mbox{ for all }\phi \in H^1(\Omega)/\R.
	$$

\begin{definition}[Weak solution to the limiting system]\label{def_weak_limit}
We shall say that the couple $(\vec{U}, \Psi)$ is a weak solution to the problem \eqref{limit1}, provided that the following properties hold:
\begin{itemize}
\item[1)]
$\vec{U} \in L^\infty(0,T;L^2 (\Omega)^3) \cap L^2(0,T;H^1(\Omega)^3)$, $\vec{U} \cdot \vec{n} |_{\partial \Omega} = \vec{0}$ for a.e. $t \in (0,T)$, $\divx \vec{U} = 0$ for a.e. $(\vec{x},t) \in \Omega \times (0,T)$ and
$\widehat\Psi \in L^{v}(0,T;Z_1) \cap H^1(0,T; M^{-1} (H^s(\Omega \times D))'),$
where $v \in [1,\infty)$ and $s > 1+ \frac{3}{2} (K +1) $, with finite relative entropy and Fisher information, i.e.,
 	\begin{equation}\label{L61_2_df}
	{\mathcal{F}}(\widehat \Psi) \in L^\infty(0,T; L^1_M(\Omega\times D)) \quad
	\mbox{ and } \quad \sqrt{\widehat\Psi} \in L^2(0,T; H^1_M(\Omega \times D)),
	\end{equation}
	\begin{equation}\label{L61_8_df}
	\tens{\tau}_1 (M \widehat\Psi)  \in L^r(\Omega\times [0,T)), \quad \mbox{ where } r \in \bigg[ 1, \frac{20}{13}  \bigg).
	\end{equation}
\item[2)]
Moreover the following relations are satisfied:
	\begin{equation}\label{weak_lim_mom}
	\begin{split}
	& - \int_0^T \int_\Omega \bar\rho \vU \cdot  \partial_t \vec{w} \dx \dt
	+ \int_0^T \int_\Omega \left[ \tens{S}(\vU) - \bar\rho \vU \otimes \vU  \right] : \nabla_x \vec{w} \dx\dt
	\\ &\qquad
	= \int_\Omega  \bar\rho \vU_0 \cdot \vec{w}(0,\cdot) \dx - \int_0^T \int_\Omega  \left( \tens{\tau}_1(M \widehat\Psi)  \right)
	: \nabla_x \vec w
	\dx\dt
	\\ &
	\hspace{-1cm}\mbox{ for all }
	 \mbox{ for all } \vec{w} \in L^{\frac{\gamma + \vartheta}{\vartheta}}(0,T; W^{1,r}(\Omega)^3),
	\mbox{ such that }  \divx \vec{w} = 0,\ \vec{w} \cdot \vec{n}|_{\partial\Omega} = 0,
\end{split}
	\end{equation}
where $\vec{U}_0  \in L^2(\Omega \times (0,T))^3$, $\vartheta(\gamma)= \frac{2\gamma - 3}{3}$ for $\frac{3}{2}<  \gamma \leq 4$
and $\vartheta(\gamma)= \frac{5}{12}\gamma$ for
$4\leq \gamma$; $r= \max \left\{  4, \frac{6\gamma}{2\gamma - 3} \right\}$, and
	\begin{equation}\label{weak_lim_FP}
	\begin{split}
	& \int_0^T \left\langle  M\partial_t \widehat\Psi , \varphi   \right\rangle_{H^s (\Omega \times D)} \dt
	+  \frac{1}{4}   \sum_{i,j=1}^K A_{ij} \int_{\Omega\times D}  M \nabla_{q_j}  {\widehat\Psi}
	\cdot  \nabla_{q_i}  {\varphi} {\rm\,d}\vec{q}\dx\dt
	\\ &\quad + \int_0^T \int_{\Omega \times D} M [ \delta \nabla_x \widehat\Psi - \vU\widehat\Psi ] \cdot \nabla_x \varphi {\rm\,d}\vec{q}\dx\dt
	- \int_0^T \int_{\Omega\times D }  M \sum_{i=1}^K [( \nabla_x \vU) \vec{q}_i ] \widehat\Psi  \cdot \nabla_{q_i} \varphi
	{\rm\,d}\vec{q}\dx\dt = 0
	\\ &\mbox{ for all } \varphi \in L^2(0,T ; H^s(\Omega \times D)),
	\end{split}
	\end{equation}
with $\widehat\Psi (\cdot, 0) = \widehat\Psi_0(\cdot)$ and $s > 1 + \frac{3}{2}(K+1).$ Here $\widehat\Psi_0 \in L^1_M(\Omega \times D;\mathbb{R}_{\geq 0})$
and $\mathcal{F}(\widehat\Psi_0) \in L^1_M(\Omega \times D;\mathbb{R}_{\geq 0})$, and $\int_D M \widehat \Psi_0 \in L^\infty(\Omega;\mathbb{R}_{\geq 0})$.
\end{itemize}
\end{definition}

Our main result is then the following theorem.

 \begin{theorem}[Main theorem]
Suppose that $\Omega \subset \R^3$ is a bounded open domain with $\partial \Omega \in C^{2, \alpha}$, with some $\alpha \in (0,1)$.
Let $p$ satisfy \eqref{p1} with $\gamma > \frac{3}{2}$, let the stress tensor $\tens{S}$ satisfy \eqref{tensS}, and let the Rouse matrix satisfy \eqref{as_AA}. Let $(\rho_\ep,\, \vec{u}_\ep,\, \psi_\ep)$ be a weak solution triple to
the system \eqref{NSFP_2} emanating from the initial data satisfying \eqref{indat_1}--\eqref{indat_3} and subject to the
boundary conditions \eqref{bc-c-slip}--\eqref{bc_psi_3}. Then, by passing to suitable subsequences if necessary, we have that
 \begin{equation}\label{conv_t_1}
 \rho_\ep \to \bar\rho \mbox{ strongly in }(L^2+L^q)(\Omega \times (0,T)) \quad \mbox{ with }q=\min\{2,\gamma\},
 \end{equation}
 \begin{equation}\label{conv_t_5}
 \vu_{0,\ep} \weak \vU_0 \mbox{ weakly in } L^2(\Omega \times (0,T))^3,
 \end{equation}
  \begin{equation}\label{conv_t_4}
 \vu_\ep \weak \vU \mbox{ weakly in } L^2(0,T; H^1(\Omega)^3),
 \end{equation}
 \begin{equation}\label{conv_t_3}
 \widehat\psi_\ep \to \widehat\Psi \mbox{ strongly in } \mbox{ in }  L^v (0,T; L^1_M(\Omega \times D)),
 \end{equation}
 \begin{equation}\label{conv_t_6}
	\varrho_\ep\weakstar \varrho
	\quad \mbox{ weakly-* in }L^\infty(0,T; L^2(\Omega )),
	\quad \varrho_\ep\weak \varrho
	\quad \mbox{ weakly in }L^2(0,T; H^1(\Omega )).
 \end{equation}
Moreover, the couple $(\vU,\Psi)$ is a weak solution to \eqref{limit1} according to Definition~\ref{def_weak_limit}
with initial datum $\vU(0,\cdot) = \vec{H}[\vU_0]$, where $\vU_0$ and $\Psi_0$ are weak limits of \eqref{indat_1} and \eqref{indat_3}, 
and $\widehat \Psi_0 :=\Psi_0/M$.
 \end{theorem}

\section{Proof of the main result}\label{sec2}
Before embarking on the proof, a remark is in order concerning our choice of boundary condition for the velocity field. 

\begin{remark}
As was noted in the Introduction, we have confined ourselves in this paper to the case of a complete slip boundary condition, which models an acoustically hard boundary. Then, if $\Omega \subset \R^3$
is a bounded domain, as is being assumed here, the gradient part of the velocity field associated to acoustic waves may exhibit fast oscillations in time,
and the convergence of the velocity field $\ue$ in the limit of $\ep \rightarrow 0$ is genuinely weak with respect to the temporal variable (see \cite{FN}). 

In the case of a Dirichlet (no-slip) boundary condition for the velocity field, a boundary layer may appear because of the presence of viscosity and the specific geometrical properties of the boundary. The decay of the acoustic waves and consequently of the velocity field can then be deduced as in \cite{Desjardins}.

\end{remark}

Our starting point is the following result concerning the existence of weak solutions to the primitive system. Its proof for the case of a homogeneous Dirichlet boundary condition for the velocity field can be found in \cite[Lemma~6.5,~6.2, Theorem~6.1]{BS2014}, and the arguments contained therein can be easily adapted to the case of complete slip boundary condition by replacing the function space of divergence-free three-component vector functions contained in $H^1_0(\Omega)$ throughout the proof in \cite{BS2014} by the function space of divergence-free three-component vector functions contained in $H^1(\Omega)$ with vanishing normal trace on $\partial\Omega$.

\begin{proposition}[Existence of solutions to the primitive system]\label{prop.eps}
For each fixed $\ep>0$ there exists a  triple $(\rho_\ep, \vu_\ep, \widehat\psi_\ep)$, which is a global weak solution to problem \eqref{NSFP_2} in the sense of Definition~\ref{weak_sol_primitiv}.
Furthermore, the solution triple $(\rho_\ep, \vu_\ep, \widehat\psi_\ep)$ satisfies, for a.e. $t' \in (0,T)$, the energy inequality
	\begin{equation}\label{L65_33}
	\begin{split}
	& \frac{1}{2}  \int_\Omega \rho_\ep(t') |\vu_\ep (t') |^2 \dx + \frac{1}{\ep^2} \int_\Omega \frac{1}{\gamma -1}p(\rho_\ep(t')) \dx
	+ k \int_{\Omega \times D} M {\mathcal{F}}(\widehat\psi_\ep(t')) {\rm\,d}\vec{q} \dx
	+ \bar\xi \| \varrho_\ep (t') \|^2_{L^2(\Omega)}
	\\ &\quad  + \mu^S c_0 \int_0^{t'} \int_0^{t'} \| \vu_\ep \|^2_{H^1(\Omega)} \dt
    \\ &\quad  + k  \int_0^{t'}  \int_{\Omega \times D} M \Big[
	\frac{a_0}{2 \lambda} \big|  \nabla_q \sqrt{\widehat\psi_\ep} \big|^2
	+ 2 \delta  \big| \nabla_x \sqrt{\widehat\psi_\ep}\big|^2   \Big] {\rm\,d}\vec{q}\dx\dt
	+2 \bar\xi \delta \int_0^{t'} \| \nabla_x \varrho_\ep \|^2_{L^2(\Omega)} \dt
	\\ & \leq
	\exp(t') \Big[
	\frac{1}{2} \int_{\Omega} \rho_{0,\ep}|\vu_{0,\ep}|^2 \dx + \int_\Omega P(\rho_{0,\ep}) \dx
	+ k  \int_{\Omega \times D} M {\mathcal{F}}(\widehat\psi_{0,\ep}) {\rm\,d}\vec{q} \dx
	+\bar\xi \int_{\Omega} \bigg(\int_D M \widehat\psi_{0,\ep} {\rm\,d}\vec{q}\bigg)^2 \dx
	\Big].
	\end{split}
	\end{equation}
 \end{proposition}

\subsection{Energy equality}

In order to obtain uniform bounds on $(\rho_\ep, \vu_\ep, \widehat\psi_\ep)$ with respect to $\ep$,  we shall prove a formal energy equality. Its derivation is in fact the same as in the proof of the existence of weak solutions in Proposition \ref{prop.eps}. For the moment we shall keep all characteristic numbers in their general form and will permit a nonzero density of body forces $\vec{f}$ so as to state the formal energy equality in its most general form. By taking the $L^2(\Omega)$ inner product of the continuity equation first with $\frac{1}{2} |\vu|^2$ and then with $P'(\rho) - P'(\bar\rho)\rho$, where $P(\rho) = \frac{p(\rho)}{\gamma - 1}$ and where $\bar\rho$ is determined by \eqref{st_st}, and taking the $L^2(\Omega)^3$ inner product of the momentum equation with $\vu$ (by noting the boundary conditions and performing partial integration) we deduce that
	\begin{equation}\label{ues_1}
	\begin{split}
	& \frac{{\rm d}}{{\rm d} t} \int_\Omega
	\bigg(\frac{1}{2} \rho |\vu |^2 + \frac{1}{{\rm Ma}^2} (P(\rho) - P'(\bar\rho)(\rho - \bar\rho) - P(\bar\rho) )
	\bigg)
	\dx + \frac{1}{{\rm Re}} \int_\Omega  \mu^S | \tens{D} \vu
	- \frac{1}{3} (\divx \vu) \tens{I} |^2 \dx
	\\ & + \frac{1}{{\rm Re}} \int_\Omega  \mu^B | \divx \vu |^2 \dx
	+ \frac{1}{{\rm Re}} \int_\Omega \tens{\tau}_1 : \tens{D} \vu \dx
	- \frac{\tilde{\xi}}{{\rm Ma}^2} \int_\Omega \varrho^2 \divx \vu \dx = \frac{1}{{\rm Fr}^2} \int_\Omega \rho \vec{f} \cdot \vu \dx .
	\end{split}
	\end{equation}
Now let us concentrate on the Fokker--Planck equation.
Multiplying the Fokker--Planck equation in \eqref{NSFP_1} by $\mathcal{F}'(\widehat\psi)$, integrating over $\Omega\times D$ and noticing that $\nabla_{q_i} M =  - M \vec{q}_i U'_i \left( \frac{ | \vec{q}_i |^2 }{2}  \right)$, we deduce that
	\begin{equation}\label{ues_2}
	\begin{split}
	& \frac{\rm d}{\dt}  \int_{\Omega \times D} M\mathcal{F}(\widehat\psi) {\rm\,d}\vec{q} \dx
	+ (K + 1) \int_{\Omega \times D}  M(\divx \vu ) \widehat\psi {\rm\,d}\vec{q}\dx\\
	& \quad- \sum_{i=1}^K \int_{\Omega\times D} (\vec{q}^T_i \vec{q}_i) M U'_i \left( \frac{ | \vec{q}_i |^2 }{2}  \right) \widehat\psi : (\nabla_x \vu) {\rm\,d}\vec{q}\dx
	\\
	& \quad
	+ 4 \delta \int_{\Omega \times D} M |\nabla_x \sqrt{\widehat\psi} |^2 {\rm\,d}\vec{q}\dx
	+ \frac{1}{{\mathrm{De}}}  \sum_{i,j=1}^K A_{ij} \int_{\Omega\times D} M \nabla_{q_i}  \sqrt{\widehat\psi}
	\cdot  \nabla_{q_j}  \sqrt{\widehat\psi} {\rm\,d}\vec{q}\dx = 0.
	\end{split}
	\end{equation}
Multiplying \eqref{ues_2} by $\frac{1-\beta}{{\rm Wi\, Re}}$  and adding to \eqref{ues_1} we deduce that
	\begin{equation}\label{ues_3}
	\begin{split}
	& \frac{{\rm d}}{{\rm d} t} \int_\Omega
	\bigg(
	 \frac{1}{2} \rho |\vu |^2 + \frac{1}{{\rm Ma}^2} (P(\rho) - P'(\bar\rho)(\rho - \bar\rho) - P(\bar\rho) )
	 \bigg)
	  \dx + \frac{1}{{\rm Re}} \int_\Omega  \mu^S | \tens{D} \vu - \frac{1}{3} (\divx \vu) \tens{I} |^2 \dx
	\\ &\quad + \frac{1}{{\rm Re}} \int_\Omega  \mu^B | \divx \vu |^2 \dx
	+ \frac{1}{{\rm Re}} \frac{1-\beta}{\mathrm{De}}  \frac{\rm d}{\dt}  \int_{\Omega \times D} M \mathcal{F}(\widehat\psi) {\rm\,d}\vec{q} \dx
	+ \frac{ 4 \delta}{\rm Re} \frac{1-\beta}{\mathrm{De}}  \delta \int_{\Omega \times D} M |\nabla_x \sqrt{\widehat\psi} |^2 {\rm\,d}\vec{q}\dx
	\\
	&\quad
	+ \frac{1}{{\rm Wi\,Re}}  \frac{1-\beta}{\mathrm{De}}  \int_{\Omega\times D}   \sum_{i,j=1}^K A_{ij} M \nabla_{q_i}  \sqrt{\widehat\psi}
	\cdot  \nabla_{q_j}  \sqrt{\widehat\psi} {\rm\,d}\vec{q}\dx
	- \frac{\tilde{\xi}}{{\rm Ma}^2} \int_\Omega \varrho^2 \divx \vu \dx
	= \frac{1}{{\rm Fr}^2} \int_\Omega \rho \vec{f} \cdot \vu 	\dx .
	\end{split}
	\end{equation}
It remains to deal with the last term on the left-hand side of \eqref{ues_3}. To this end we  integrate the Fokker--Planck equation over $D$ and we get \eqref{rem_11}. After multiplying \eqref{rem_11} by $\varrho$ and integrating over $\Omega$  we have
	\begin{equation}\label{ues_4}
	\frac{{\rm d}}{{\rm d} t} \int_\Omega \varrho^2 \dx + 2 \delta \int_\Omega |\nabla_x \varrho|^2 \dx
	= -\int_\Omega \varrho^2 (\divx \vu ) \dx .
	\end{equation}
Multiplying \eqref{ues_4} by $\frac{\tilde\xi}{{\rm Ma}^2}$ and substituting the resulting expression into \eqref{ues_3}  gives
	\begin{equation}\label{ues_5}
	\begin{split}
	& \frac{\rm d}{\dt} \left[ \int_\Omega
	\left(
	 \frac{1}{2} \rho |\vu|^2 + \frac{1}{{\rm Ma}^2}  (P(\rho) - P'(\bar\rho)(\rho - \bar\rho) - P(\bar\rho) )
	 \right)
	 \dx
	+ \frac{1-\beta}{{\rm Re}\,{\mathrm{De}}} \int_{\Omega \times D} M \mathcal{F}(\widehat\psi){\rm\,d}\vec{q}\dx+ \frac{\tilde\xi}{{\rm Ma}^2} \int_\Omega  \varrho^2 \dx
	 \right]
	\\ &\quad
	+ \Big\{
	 \frac{1}{{\rm Re}} \int_\Omega \mu^S | \tD (\vu ) - \frac{1}{3} ( \divx \vu ) |^2 \dx
	 + \frac{1}{{\rm Re} } \int_\Omega \mu^B | \divx \vu |^2 \dx
	 + \frac{4\delta(1-\beta)}{{\rm Re}\, {\mathrm{De}}} \int_{\Omega \times D} M | \nabla_x \sqrt{\widehat\psi}|^2 {\rm\,d}\vec{q}\dx
	 \\ &\quad
	+ \frac{1-\beta}{{\rm Re}\,  {\mathrm{De}}^2}
	 \int_{\Omega \times D} M \sum_{i,j = 1}^K A_{ij} \nabla_{q_i}  \sqrt{\widehat\psi} \cdot  \nabla_{q_j}  \sqrt{\widehat\psi} {\rm\,d}\vec{q}\dx
	 + \frac{2 \delta \tilde\xi}{{\rm Ma}^2} \int _\Omega | \nabla_x \varrho |^2 \dx
	 \Big\}
	 \\ &
	 = \frac{1}{{\rm Fr}^2} \int_\Omega  \rho \vec{f} \cdot \vu \dx .
	\end{split}
	\end{equation}

In accordance with our assumptions on the characteristic numbers we shall take $\vec{f} = \vec{0}$, ${\rm Ma} = \ep$, ${\mathrm{De}} = 1$, ${\rm Re} = 1$, and $\tilde\xi = \bar\xi \ep^2$ in the above equality; thereby,
	\begin{equation}\label{ues_6}
	\begin{split}
	& \frac{{\rm d}}{{\rm d} t} \left[ \int_\Omega
	\left(
	\frac{1}{2} \rho_\ep |\vu_\ep|^2 + \frac{1}{{\ep}^2}  (P(\rho_\ep) - P'(\bar\rho)(\rho_\ep - \bar\rho) - P(\bar\rho) )
	\right)
	 \dx
	+ (1-\beta) \int_{\Omega \times D} M \mathcal{F}(\widehat\psi_\ep){\rm\,d}\vec{q}\dx+ {\bar\xi} \int_\Omega  \varrho^2_\ep \dx
	 \right]
	\\ &\quad
	+ \Big\{
	 \int_\Omega \mu^S | \tD (\vu_\ep ) - \frac{1}{3} ( \divx \vu_\ep ) |^2 \dx
	 + \int_\Omega \mu^B | \divx \vu_\ep |^2 \dx
	 + 4\delta(1-\beta) \int_{\Omega \times D} M | \nabla_x \sqrt{\widehat\psi_\ep}|^2 {\rm\,d}\vec{q}\dx
	 \\ & \quad
	+ (1-\beta)
	 \int_{\Omega \times D} M \sum_{i,j = 1}^K A_{ij} \nabla_{q_i}  \sqrt{\widehat\psi_\ep} \cdot  \nabla_{q_j}  \sqrt{\widehat\psi_\ep} {\rm\,d}\vec{q}\dx
	 + 2 \delta \bar\xi \int _\Omega | \nabla_x \varrho_\ep |^2 \dx
	 \Big\}
	 \\ &
	 = 0 \, .
	\end{split}
	\end{equation}
Consequently, using also \eqref{as_AA} we obtain
		\begin{equation}\label{ues_7}
	\begin{split}
	&\left[ \int_\Omega
	\left(
	\frac{1}{2} \rho_\ep |\vu_\ep|^2  + \frac{1}{{\ep}^2}  ({P}(\rho_\ep)
	- {P}'(\bar\rho)(\rho_\ep - \bar\rho) - {P}(\bar\rho) )
	\right)
	 \dx
	+ (1-\beta) \int_{\Omega \times D} M \mathcal{F}(\widehat\psi_\ep){\rm\,d}\vec{q}\dx+ {\bar\xi} \int_\Omega  \varrho^2_\ep \dx
	 \right] (t)
	\\ &\quad
	+ \int_0^t  \Big\{
	\int_\Omega \mu^S | \tD (\vu_\ep ) - \frac{1}{3} ( \divx \vu_\ep ) |^2 \dx
	 + \int_\Omega \mu^B | \divx \vu_\ep |^2 \dx
	 + 4\delta(1-\beta) \int_{\Omega \times D} M | \nabla_x \sqrt{\widehat\psi_\ep}|^2 {\rm\,d}\vec{q}\dx
	 \\ &\quad
	+ (1-\beta) a_0
	 \int_{\Omega \times D} M| \nabla_q \sqrt{\widehat\psi_\ep}|^2 {\rm\,d}\vec{q}\dx
	 + 2 \delta \bar\xi \int _\Omega | \nabla_x \varrho_\ep |^2 \dx
	 \Big\} \dt
	 \\ &
	 \leq \left[ \int_\Omega  \Big(\frac{1}{2} \rho_{0,\ep} |\vu_{0,\ep}|^2  + \frac{1}{{\ep}^2}  ({P}(\rho_{0,\ep}) - {P}'(\bar\rho)(\rho_{0,\ep} - \bar\rho) - {P}(\bar\rho) ) \Big) \dx \right.
    \\ & \quad \left.
 	+ (1-\beta) \int_{\Omega \times D} M \mathcal{F}(\widehat\psi_{0,\ep}){\rm\,d}\vec{q}\dx+ {\bar\xi} \int_\Omega  \varrho^2_{0,\ep} \dx
	 \right], \qquad \mbox{for a.e. $t \in (0,T)$}.
	\end{split}
	\end{equation}

\begin{remark}\label{cov_1}	
We note that ${P}''(\rho) =  p'(\rho) / \rho > 0$ for all $\rho \in (0,\infty)$. Thus ${P}$ is a strictly convex function of $\rho$ on $[0,\infty)$, and
${P}(\rho) - {P}'(\rs)(\rho - \rs) - {P}(\rs) \approx c (\rho - \bar\rho)^2$ provided that $\rho$ is close to $\bar\rho$.
\end{remark}

\subsection{Uniform estimates}\label{ss:uni}
In order to deduce uniform bounds from \eqref{ues_6} we introduce, similarly as in \cite{FN}, the decomposition into essential and residual parts of a measurable function $h$ as
	\begin{equation*}\label{ess-def}
	h = [h]_{\ess} + [h]_{\res},\quad\  [ h]_{\ess} := \chi(\rho_\ep )h,\quad \  [h]_{\res} := (1-\chi(\rho_\ep))h,
	\end{equation*}
with $\chi \in C^\infty_c ( (0,\infty) )$,
	$0\leq \chi \leq 1, \ \chi=1$  on the set $ {\cal{O}}_{\ess}$,
where
	\begin{equation*}\label{ess-set}
	 {\cal{O}}_{\ess} :=  [\vrs/2 , 2 \vrs ] ,\quad
	 {\cal{O}}_{\res} := (0,\infty)\setminus {\cal{O}}_{\ess}.
 	\end{equation*}

The following estimates result from the energy estimate, with a constant $c$ independent of $\ep$.
According to \eqref{indat_1}--\eqref{indat_3} the expression on the right-hand side
of \eqref{ues_7} is bounded for $\ep \to 0$.
Thus we infer that
	\begin{equation}\label{ue_1}
	\esssup_{t \in (0,T)} \int_\Omega \re | \ue |^2 \dx \leq c ,
	\end{equation}
	\begin{equation}\label{ue_2}
	\int_0^T \int_\Omega \mu^S | \tD (\ue ) - \frac{1}{3} ( \divx \ue ) |^2 \dx
	 + \int_\Omega \mu^B | \divx \ue |^2 \dx \leq c ,
	 \end{equation}
	 \begin{equation}\label{ue_3}
	 \esssup_{t \in (0,T)} \int_\Omega \left[ (\tilde{P}(\re) - \tilde{P}'(\bar\rho)(\re - \bar\rho) - \tilde{P}(\bar\rho) ) \right]_{\rm ess} \dx \leq c\,\ep^2,
	 \end{equation}
		 \begin{equation}\label{ue_4}
	 \esssup_{t \in (0,T)} \int_\Omega \left[ (\tilde{P}(\re) - \tilde{P}'(\bar\re)(\rho - \bar\rho) - \tilde{P}(\bar\rho) ) \right]_{\rm res} \dx \leq c\,\ep^2,
	 \end{equation}	
	\begin{equation}\label{ue_42}
	\esssup_{t \in (0,T)} \int_\Omega | \varrho_\ep|^2 \dx + \int_0^T \int_\Omega | \nabla_x \varrho_\ep|^2 \dx\dt \leq c,
	\end{equation}
	\begin{equation}\label{ub_p_2}
	\| {\mathcal{F}}(\widehat\psi_\ep )\|_{L^\infty(0,T; L^1_M(\Omega\times D))} \leq c,
	\end{equation}
	\begin{equation}\label{ub_p_2_0}
	\| \nabla_x \sqrt{\widehat\psi_\ep} \|_{L^2(0,T; L_M^2(\Omega\times D))}
	+ \| \nabla_q \sqrt{\widehat\psi_\ep } \|_{L^2(0,T; L_M^2(\Omega\times D))} \leq c.
	\end{equation}
By \eqref{ub_p_2_0}, together with \eqref{ue_42}, we have that
	\begin{equation}\label{ub_p_2_2}
	\| \sqrt{\widehat\psi_\ep } \|_{L^2(0,T; H^1_M(\Omega \times D))} \leq c .
	\end{equation}
The bound \eqref{ue_42} and Sobolev embedding yield
	\begin{equation}\label{ue_10}
	 \| \varrho_\ep \|_{L^2(0,T; L^6(\Omega))} \leq c.
	\end{equation}
Hence, by an interpolation argument, we deduce that
	\begin{equation}\label{ue_10_2}
	 \| \varrho_\ep \|_{L^{\frac{10}{3}}((0,T) \times \Omega))}
	 + \| \varrho_\ep \|_{L^4(0,T; L^3(\Omega))} \leq c .
	\end{equation}
As $\tilde{P}$ is strictly convex, the relation \eqref{ue_3} yields
	\begin{equation}\label{ue_5}
	\esssup_{t \in (0,T)} \int_\Omega \left| \left[  \frac{\re - \rs}{\ep} \right]_{\rm ess}\right|^2 \dx \leq  c,
	\end{equation}
and, by virtue of \eqref{p1} and \eqref{ue_4},	
	\begin{equation}\label{ue_6}
	\esssup_{t \in (0,T)} \int_\Omega [1+  \re ]^\gamma_{\rm res} \dx \leq c\, \ep^2.
	\end{equation}
Next, as a direct consequence of \eqref{ue_1}, we have that
	\begin{equation}\label{ue_7}
	\esssup_{t \in (0,T)} \| [\rho_\ep \ue ]_{\rm ess} \|_{L^2(\Omega)} \leq c,
	\end{equation}
and by Korn's inequality  (\cite[Theorem~10.17]{FN}) and  \eqref{ue_7}, \eqref{ue_2} gives
	\begin{equation}\label{ue_8}
	\| \ue \|_{L^2(0,T; H^{1}(\Omega))} \leq c.
	\end{equation}
Moreover by \eqref{ue_1}, \eqref{ue_3}, \eqref{ue_4}, and \eqref{ue_8} we have that
\begin{equation}\label{ub_p_3}
	\begin{split}
	& \| \rho_\ep \|_{L^\infty (0,T; L^\gamma(\Omega))}
	+ \| \rho_\ep \vu_\ep \|_{L^\infty(0,T; L^{ \frac{2\gamma}{\gamma + 1}}(\Omega))}
	+  \| \rho_\ep \vu_\ep \|_{L^2(0,T; L^{ \frac{6\gamma}{\gamma + 6}}(\Omega))}
	+  \| \rho_\ep |\vu_\ep|^2 \|_{L^2(0,T; L^{ \frac{6\gamma}{4\gamma + 3}}(\Omega))}
	\leq c.
	\end{split}
	\end{equation}
Furthermore, by using the continuity equation we additionally infer that
		\begin{equation*}
	\| \partial_t \rho_\ep\|_{L^2(0,T; W^{1,s'} (\Omega)')} \leq c,\quad \mbox{ where  $s$ is as in \eqref{T61_1}}.
	\end{equation*}

Next, we shall establish the necessary bounds on the extra stress tensor $\tens{\tau}_1$.
We deduce from \eqref{C1}, \eqref{Max1}, and by noting that $M |_{\partial D} = 0$, that
	\begin{equation}\label{p422}
	\tens{C}_i (M\varphi) = -\int_{D} (\nabla_{q_i} M) \vec{q}_i^T \varphi {\rm\,d}\vec{q}
	= \int_{D}  M (\nabla_{q_i} \varphi)  \vec{q}_i^T {\rm\,d}\vec{q} + \left(  \int_D  M \varphi {\rm\,d}\vec{q} \right) \tens{I} .
	\end{equation}
As $\nabla_{{q_i}} \varphi = \nabla_{q_i} (\sqrt{\varphi}\,)^2 = 2 \sqrt{\varphi}\, \nabla_{q_i} \sqrt{\varphi}$ for any {sufficiently smooth} nonnegative function $\varphi$, we have that
	\begin{equation}\label{p423}
	\begin{split}
	\|\tens{C}_i (M\varphi) \|_{L^r(\Omega)}
	& \leq c \left[
	\int_\Omega \left( \int_D M \varphi {\rm\,d}\vec{q}  \right)^{\frac{r}{2}}
	\left( \int_D M | \nabla_{q_i} \sqrt{\varphi}|^2 {\rm\,d}\vec{q}  \right)^{\frac{r}{2}} \dx
	+ \int_\Omega \left( \int_D M \varphi {\rm\,d}\vec{q}   \right)^r
	\right]^{\frac{1}{r}}\\
	&
	\leq c \left[
	\| \nabla_{q_i} \sqrt{\varphi} \|_{L^2_M(\Omega \times D)}
	\left\| \int_D M \varphi {\rm\,d}\vec{q}  \right\|^{\frac{1}{2}}_{L^{\frac{r}{2-r}}(\Omega)}
	+ \left\| \int_D M \varphi {\rm\,d}\vec{q}  \right\|_{L^r(\Omega)}
	\right],
	\end{split}
	\end{equation}
for $r\in[1,2)$. Then, for $s\in [1,2]$ we get that, for any such function $\varphi$,
	\begin{equation}\label{p424}
	\| \tens{C}_i(M \varphi) \|_{L^s(0,T; L^r(\Omega))}
	\!\leq\!
	c\left[
	\| \nabla_{q_i} \sqrt{\varphi} \|_{L^2(0,T;L^2_M(\Omega\times D)}
	\left\|  \int_D M \varphi {\rm\,d}\vec{q}  \right\|^{\frac{1}{2}}_{L^v(0.T;L^{\frac{r}{2-r}}(\Omega))}
	+ \left\|  \int_D M \varphi {\rm\,d}\vec{q}  \right\|_{L^s(0.T;L^r(\Omega))}
 	\right],
	\end{equation}
where $v= \frac{s}{2-s}$ if $s\in [1,2)$ and $v = \infty $ if $s=2$. We deduce from \eqref{ue_42}, \eqref{p423} and \eqref{ub_p_2_0}
that, for $i=1, \dots, K,$
	\begin{equation}\label{p425}
	\| \tens{C}_i (M \widehat{\psi}_\ep) \|_{L^s(0,T;L^r(\Omega))} \leq c, \quad
	\mbox{ as } \quad \| \varrho_\ep \|_{L^v(0,T;L^{\frac{r}{2-r}}(\Omega))} \leq c,
	\end{equation}
where $r\in [1,2),$ $s\in [1,2]$, and $v= \frac{s}{2-s}$ if $s \in [1,2)$, and $v = \infty$ if $s=2$.
Next by \eqref{ue_42}  and interpolation we have
	\begin{equation}\label{p426}
	\| \varrho_\ep \|_{L^{\frac{2}{v}} (0,T;L^v(\Omega))}
	\leq c \|\varrho_\ep \|^{1-\nu}_{L^\infty(0,T;L^2(\Omega))} \| \varrho \|^\nu_{L^2(0,T;H^1(\Omega))} \leq c,
	\end{equation}
where $\nu = \frac{3(v-2)}{2v}$, and $v\in (2,6]$. By
\eqref{def_t1}, \eqref{ue_42},  \eqref{ub_p_2_0},
\eqref{ue_10}, \eqref{ue_10_2}
\eqref{p424},
\eqref{p428a}	we deduce that
	\begin{equation}\label{p428a}
	\| \tens{C}_i (M \widehat{\psi}_\ep ) \|_{L^2(0,T;L^{\frac{4}{3}} (\Omega))}
	+ \| \tens{C}_i (M \widehat{\psi}_\ep ) \|_{L^{\frac{20}{13}} (\Omega \times (0,T) )}  \leq c,
	\end{equation}
	\begin{equation}\label{p428b}
	\| \tens{\tau}_1 (M \widehat{\psi}_\ep ) \|_{L^2(0,T;L^{\frac{4}{3}} (\Omega))}
	+ \| \tens{\tau}_1 (M \widehat{\psi}_\ep ) \|_{L^{\frac{20}{13}} (\Omega \times (0,T) )}
	+ \| \vec{\tau}_1(M\widehat\psi_\ep ) \|_{L^{\frac{4}{3}}(0,T; L^{\frac{12}{7}}(\Omega))}
	 \leq c ,
	\end{equation}
	where $c$ is independent of $\ep$.

Now let us show that
	\begin{equation}\label{lem_46}
	\| M \partial_t \widehat{\psi}_\ep \|_{L^2(0,T; H^s(\Omega \times D)')} \leq c,
	\end{equation}
with $s > 1+ \frac{3}{2}(K+1) $.
Testing the Fokker--Planck equation  \eqref{NSFP_2}$_3$ with $\varphi \in L^2 (0,T; W^{1,\infty}(\Omega \times D))$ and since
$\nabla_x \psi = 2 \sqrt{\psi}\, \nabla_x \sqrt{\psi}$ for sufficiently smooth $\psi$, we infer  that
	\begin{equation}\label{pt_psi}
	\begin{split}
	\left|
	\int_0^T \int_{\Omega \times D} M \partial_t \widehat{\psi}_\ep \varphi {\rm\,d}\vec{q}\dx \dt
	\right|
	& \leq
	2 \delta \left|
	\int_0^T\int_{\Omega \times D} M \sqrt{\widehat{\psi}_\ep} \nabla_x \sqrt{\widehat{\psi}_\ep}
	\cdot \nabla_x \varphi {\rm\,d}\vec{q} \dx \dt
	\right|
	\\ & \quad +
	\frac{1}{2} \left|
	\sum_{i,j =1}^K A_{i,j} \int_0^T \int_{\Omega\times D}
	M \sqrt{\widehat{\psi}_\ep} \nabla_{q_j} \sqrt{\widehat{\psi}_\ep}
	\cdot \nabla_{q_i} \varphi {\rm\,d}\vec{q}\dx\dt
	\right|
	\\ & \quad
	+ \left|
	\int_0^T \int_{\Omega \times D} M \vu_\ep \widehat{\psi}_\ep \cdot  \nabla_x \varphi {\rm\,d}\vec{q}\dx\dt
	\right|
	\\ & \quad
	+
	 \left|
	\int_0^T \int_{\Omega \times D}
	M  \sum_{i=1}^K \left[ (\nabla_x \vu_\ep) \vec{q}_i \right] \widehat{\psi}_\ep  \cdot \nabla_{q_i} \varphi {\rm\,d}\vec{q}\dx\dt
	\right|
	\\ &
	\leq c \| \varrho_\ep \|_{L^\infty(0,T;L^2(\Omega))}
	\Big[
	\| \nabla_x \sqrt{{\widehat\psi}_\ep} \|_{L^2(0,T; L^2_M(\Omega \times D))}
	+ \| \nabla_q \sqrt{{\widehat\psi}_\ep} \|_{L^2(L_M^2 (\Omega \times D))}
	\\ & \quad + \| \nabla_x \vu_\ep \|_{L^2(0,T;L^2(\Omega))} + \|  \vu_\ep \|_{L^2(0,T;L^2(\Omega))}
	\Big] \| \varphi\|_{L^2(0,T; W^{1,\infty}(\Omega \times D))}
	\\ & \leq c.
	\end{split}
	\end{equation}
The last inequality holds by  \eqref{ue_42}, \eqref{ub_p_2_0}, \eqref{ue_8} and as
$H^s (\Omega \times D) \subset W^{1,\infty} (\Omega \times D)$ thanks to our assumption that $s>
1 + \frac{3}{2}(K+1)$.
Therefore \eqref{lem_46} holds.

Finally, we note that by choosing $\varphi = 1 $ in \eqref{T61_3} and by the choice of the initial datum $\widehat\psi_{0,\ep}$ (cf. \eqref{indat_3}),
we have that
	\begin{equation}\label{429}
	\int_\Omega \varrho_\ep \dx = \int_{\Omega \times D} M \widehat\psi_{\ep} \dq \dx = \int_{\Omega \times D} M \widehat\psi_{0,\ep} \dq \dx \leq c,
	\end{equation}
where again $c$ is independent of $\ep$.

\subsection{Convergence as $\ep \rightarrow 0$}
Obviously $\rho_\ep - \bar\rho = [\rho_\ep - \bar\rho]_{ess} + [\rho_\ep - \bar\rho]_{res}$; hence, thanks to \eqref{ue_5}, \eqref{ue_6}, we have that
	\begin{equation}\label{conv_1}
	\rho_\ep \to \bar\rho \quad \mbox{ strongly in } L^\infty(0,T;(L^\gamma + L^2)(\Omega)) .
	\end{equation}
Next, by \eqref{ue_8} we have
	\begin{equation}\label{conv_2}
	\vu_\ep \weak \vU \quad \mbox{ weakly in } L^2(0,T;H^{1}(\Omega)^3)
	\end{equation}
(by extracting a subsequence if necessarily).
Using \eqref{ue_1}, \eqref{ue_7} and Sobolev embedding we get
	\begin{equation}\label{conv_3_0}
	\rho_\ep \vu_\ep \weakstar \bar\rho\, \vU \quad \mbox{ weakly-* in } L^\infty(0,T;L^{\frac{2\gamma}{1+ \gamma}}(\Omega)^3),
	\end{equation}
	\begin{equation}\label{conv_3}
	\rho_\ep \vu_\ep \weak \bar\rho\, \vU \quad \mbox{ weakly in } L^2(0,T;L^{\frac{6\gamma}{6+ \gamma}}(\Omega)^3) .
	\end{equation}
Moreover, from  \eqref{ue_1} and \eqref{ue_6} it follows that
	\begin{equation}\label{ue_9}
	[\re \ue ]_{\rm res } \to 0 \quad \mbox{ in } L^\infty(0,T; L^s (\Omega)^3)\quad  \mbox{ as } \ep \to 0 \mbox{ for } 1\leq s \leq 2\gamma / (\gamma + 1).
	\end{equation}
Then, employing \eqref{conv_1} and \eqref{conv_2}, we have from \eqref{T61_1}  that
$$\int_0^T\int_\Omega \vU \cdot \nabla_x \eta \dx\, \mathrm{d}t = 0 \quad \mbox{ for all } \eta \in C_c^\infty(\overline\Omega \times (0,T)). $$
Since the limiting velocity $\vU \in L^2(0,T; H^1(\Omega)^3)$, we obtain that
	\begin{equation}\label{divU}
	\divx \vU = 0 \quad \mbox{ for a.e. } (\vec{x},t) \in \Omega \times (0,T) \mbox{ and } \vU \cdot \vec{n}|_{\partial\Omega } = 0 \mbox{ in the sense of traces}.
	\end{equation}
Thus in the limit of $\ep \to 0$ the continuity equation is reduced to the divergence-free condition for the limiting velocity field $\vU$. This justifies the use of solenoidal test functions for the momentum equation.

As a direct consequence of \eqref{conv_2} and \eqref{divU} we have
	\begin{equation}\label{limit_S}
	\tS ( \vu_\ep ) \weak \mu^S \tens{D} \vU \quad \mbox{ weakly in } L^2(\Omega \times (0,T))^{3 \times 3}.
	\end{equation}
Thanks to \eqref{conv_2} and \eqref{conv_3_0} we deduce that
	\begin{equation}\label{convectiv_01}
	\rho_\ep \vu_\ep \otimes \vu_\ep \weak \overline{\bar\rho \vU \otimes \vU}
	\quad \mbox{ weakly in } L^2(0,T; L^{\frac{6\gamma}{4\gamma + 3}} (\Omega)^{3\times 3}).
	\end{equation}
In the above ${\overline{\overline{\rho} \vU \otimes \vU }}$ denotes a weak limit of $\{ \rho_\ep \vu_\ep \otimes \vu_\ep \}_{\ep >0}$ in $L^2(0,T;L^q(\Omega)^{3\times 3})$ for a certain $q>1$. Although we do not really expect that
	$${\overline{\bar\rho \vU \otimes \vU }} = \bar\rho \vU \otimes \vU,$$
we will show in the next section that
	\begin{equation}\label{co_1}
	\int_0^T \int_\Omega {\overline{\bar\rho \vU \otimes \vU }} : \nabla_x \vec{w} \dx\dt =
	\int_0^T \int_\Omega {[{\bar\rho \vU \otimes \vU }]} : \nabla_x \vec{w} \dx\dt
	\end{equation}	
for any $\vec{w} \in C^\infty_c((0,T) \times \overline{\Omega})^3$, $\divx \vec{w} = 0$, $\vec{w}\cdot \vec{n}|_{\partial\Omega} = 0$.
The relation \eqref{co_1} may be interpreted as an expression of the fact that the  difference $\divx(\overline{\bar\rho \vU \otimes \vU } - {\bar\rho \vU \otimes \vU })$ is proportional to a gradient and that it can be therefore incorporated into the limiting pressure.

\bigskip

Next we turn our attention to passage to the limit $\ep \rightarrow 0$ in the terms related to the probability density function $\widehat\psi_\ep$. In particular, we will show the strong convergence of the sequence $\{ \widehat\psi_\ep \}_{\ep>0}$
using Dubinski\v\i's compactness theorem (cf. \cite{BS2012}, for example,) stated in the next lemma.
\begin{lemma}\label{Dub}
Let $X_1$ be a seminormed set, and let $X_2$ and $X_3$ be  Banach spaces such that $X_1 \subset\subset X_2 \subset X_3$. Then, for $p_1,\,p_2 \in [1,\infty) $ the following embedding is compact
	$$
	\left\{f\in L^{p_1}(0,T;X_1) \ :  \partial_t f \in L^{p_2}(0,T;X_3) \right\} \subset\subset
	L^{p_1} (0,T; X_2).
	$$
\end{lemma}

We shall employ Lemma \ref{Dub} with $X_2 = L^1_M(\Omega \times D)$, $X_3 = M^{-1}H^s(\Omega \times D)'$,
	$$X_1 = \left\{ \phi \in Z_1 \ : \
	\int_{\Omega\times D}  M \left[ | \nabla_q \sqrt{\phi}|^2 + | \nabla_x \sqrt{\phi} |^2
	\right] {\rm\,d}\vec{q}\dx < \infty
	\right\}$$
$p_1 = v$ and $p_2 = 2$. We note that $X_1 \subset\subset X_2$ (for the details, see  \cite[Section~5]{BS10}). Moreover, $X_2 \subset X_3$. Hence,  thanks to \eqref{lem_46} and \eqref{ub_p_2_2}, we have
	\begin{equation}\label{cc4}
	 \widehat\psi_\ep \to \widehat\Psi \quad \mbox{ in }  L^v (0,T; L^1_M(\Omega \times D)).
	 \end{equation}
We infer also that
 	\begin{equation}\label{L61_4}
	M^{\frac{1}{2}} \nabla_x \sqrt{\widehat\psi_\ep} \weak M^{\frac{1}{2}} \nabla_x \sqrt{\widehat\Psi}
	\quad \mbox{ weakly in }L^2(0,T; L^2(\Omega \times D )^3),
	\end{equation}
	\begin{equation}\label{L61_5}
	M^{\frac{1}{2}} \nabla_q \sqrt{\widehat\psi_\ep} \weak M^{\frac{1}{2}} \nabla_q \sqrt{\widehat\Psi}
	\quad \mbox{ weakly in }L^2(0,T; L^2(\Omega \times D )^{3K}),
	\end{equation}
	\begin{equation}\label{L61_6}
	M \frac{\partial \widehat\psi_\ep}{\partial t} \weak M \frac{\partial \widehat\Psi}{\partial t}
	\quad \mbox{ weakly in }L^2(0,T; H^s(\Omega \times  D )').
	\end{equation}
By \eqref{cc4}, \eqref{ue_42}, and \eqref{p428b} we obtain that
	\begin{equation}\label{L61_8_c}
	\tens{\tau} (M \widehat\psi_\ep) \to \tens{\tau} (M \widehat\Psi)
	\quad \mbox{ strongly in }L^r(\Omega \times (0,T))^{3 \times 3},
	\quad \mbox{ where }r \in \bigg[ 1, \frac{20}{13}  \bigg).
	\end{equation}
We observe that, by \eqref{ue_42}, as $\ep \to 0$ we have that
	\begin{equation}\label{L61_11}
	\varrho_\ep\weakstar \varrho
	\quad \mbox{ weakly-* in }L^\infty(0,T; L^2(\Omega )),
	\quad \varrho_\ep\weak \varrho
	\quad \mbox{ weakly in }L^2(0,T; H^1(\Omega )).
	\end{equation}
Next, we note that
	\begin{equation}\label{L61_10_0}
	\varrho = \int_D M\widehat\Psi {\rm\,d}\vec{q} \in L^\infty(0,T;L^2(\Omega)) \cap L^2(0,T; H^1(\Omega)).
	\end{equation}
Indeed, \eqref{L61_10} and \eqref{cc4}, together with \eqref{L61_11} provide \eqref{L61_10_0}.

By \eqref{conv_2}, \eqref{conv_3}, \eqref{L61_8_c} and \eqref{L61_10}, we may pass to the limit $\ep \to 0$
in \eqref{T61_2} with test functions $\vec{w} \in C^1_c([0,T);C^\infty (\overline\Omega)^3)$, $\vec{w}\cdot \vec{n}|_{\partial\Omega}=0$ and $\divx\vec{w}=0$, and obtain that
	\begin{equation}\label{weak_lim_mom_0}
	\begin{split}
	& \int_0^T  \int_\Omega \bar\rho \vU \partial_t \vec{w} \dx\dt
	+ \int_0^T \int_\Omega  \left[  \overline{\bar\rho \vU \otimes \vU }
	- \mu^S \tens{D} \vU  \right] : \nabla_x \vec{w}  \dx\dt
	\\ &=  \int_0^T \int_\Omega  \left( \tens{\tau}_1(M \widehat\Psi)  \right)
	: \nabla_x \vec w
	\dx\dt
	- \int_\Omega  \bar\rho \vU_0 \vec{w}(0) \dx\dt,
	\end{split}
	\end{equation}
where we have assumed that $\vu_{\ep,0} \weak \vU_0$ in $L^2(\Omega)^3$. We note that, since in the above formulation \eqref{weak_lim_mom_0} we have restricted ourselves to divergence-free test functions, the limiting term $\nabla_x \varrho^2$ does not appear because the gradient of a scalar function can be absorbed into the limiting pressure $\Pi$.

Next let us  pass to the limit in \eqref{T61_3} and finally  obtain \eqref{weak_lim_FP}.
Initially we fix the test function to $\phi \in C([0,T]; C^\infty (\overline{\Omega \times D}))$. The first term of
\eqref{T61_3} converges to the first term of \eqref{weak_lim_FP} thanks to \eqref{L61_6}.
For the second term of \eqref{T61_3} we observe that
	\begin{equation*}
	\begin{split}
	\int_0^T \int_{\Omega \times D} M \nabla_{q_i} \widehat{\psi}_\ep \cdot 	\nabla_{q_j} \phi {\rm\,d}\vec{q}\dx\dt
	 =  &~ 2 \int_0^T \int_{\Omega \times D} M \left( \sqrt{\widehat{\psi}_\ep} - \sqrt{\widehat{\Psi}} \right)
	\nabla_{q_i} \sqrt{\widehat{\psi}_\ep} \cdot \nabla_{q_j} \phi {\rm\,d}\vec{q}\dx\dt \\
	&+ 2 \int_0^T \int_{\Omega \times D} M \sqrt{\widehat{\Psi}}\,
	\nabla_{q_i} \sqrt{\widehat{\psi}_\ep} \cdot \nabla_{q_j} \phi {\rm\,d}\vec{q}\dx\dt = : I_1 + I_2.
	 \end{split}
	\end{equation*}
Hence, by noting that $| \sqrt{c_1} - \sqrt{c_2} | \leq \sqrt{|c_1 - c_2|}$ for all $c_1,\,c_2 \in \R_{\geq 0}$ and by
\eqref{ub_p_2_2}, we have that
	\begin{equation*}
	\begin{split}
	|I_1|  & \leq C \left\| \sqrt{\widehat{\psi}_\ep} - \sqrt{\widehat{\Psi}}  \right\|_{L^2(0,T;L^2_M(\Omega \times D ))}
	\| \nabla_{q_j} \phi \|_{L^\infty(0,T;L^\infty(\Omega \times D))}\\
	& \leq C \left\| {\widehat{\psi}_\ep} - {\widehat{\Psi}}  \right\|_{L^1(0,T;L^1_M(\Omega \times D ))}
	\| \nabla_{q_j} \phi \|_{L^\infty(0,T;L^\infty(\Omega \times D))}.
	\end{split}
	\end{equation*}
Thus, \eqref{cc4} implies that $I_1$ converges to zero as $\ep \to 0$. Similarly, because
$M^{\frac{1}{2}} \sqrt{\widehat\Psi}\, \nabla_{q_j} \phi  \in L^2(0,T;L^2(\Omega\times D)^3)$ for $j=1,\dots,K$,
it follows form \eqref{L61_5} that, as $\ep \to 0$,
	\begin{equation*}
	I_2 \to 2 \int_0^T \int_{\Omega \times D} M  \sqrt{\widehat{\Psi}}\,
	\nabla_{q_i} \sqrt{\widehat{\Psi}} \cdot \nabla_{q_j} \phi {\rm\,d}\vec{q}\dx\dt
	=  \int_0^T \int_{\Omega \times D} M
	\nabla_{q_i} {\widehat{\Psi}} \cdot \nabla_{q_j} \phi {\rm\,d}\vec{q}\dx\dt.
	\end{equation*}
Therefore the second term in \eqref{T61_3} converges to the second term in  \eqref{weak_lim_FP}.
For the remaining terms in  \eqref{T61_3} we note that
	\begin{equation*}
	\begin{split}
	\int_0^T \int_{\Omega \times D}
	[M (\nabla_x \vu_\ep ) \vec{q}_i ] \widehat\psi_\ep \cdot \nabla_{q_j} \phi {\rm\,d}\vec{q}\dx\dt
	 = & \int_0^T \int_{\Omega \times D}
	[M (\nabla_x \vu_\ep ) \vec{q}_i ] \left( \widehat\psi_\ep - \widehat\Psi \right) \cdot \nabla_{q_j} \phi {\rm\,d}\vec{q}\dx\dt
	\\ & + \int_0^T \int_{\Omega \times D}
	[M (\nabla_x \vu_\ep ) \vec{q}_i ] \widehat\Psi \cdot \nabla_{q_j} \phi {\rm\,d}\vec{q}\dx\dt
	\\
	 = : & ~I_3 + I_4.
	\end{split}
	\end{equation*}
Next, by \eqref{ue_8}, using Sobolev embedding and  \eqref{429}, and recalling \eqref{L61_10}, we obtain that
	\begin{equation*}
	\begin{split}
	|I_3| & \leq c\left\|   \int_{ D}
	M |  \widehat\psi_\ep - \widehat{\Psi} | {\rm\,d}\vec{q}  \right\|_{L^2(0,T;L^2(\Omega))}
	  \| \nabla_{q_j} \phi \|_{L^\infty(0,T;L^\infty(\Omega \times D))}	   \\
	  & \leq c
	\left\|    \widehat\psi_\ep - \widehat{\Psi} \right\|^{\frac{2}{5}}_{L^2(0,T;L^1_M(\Omega \times D))}
	\| \varrho_\ep + \varrho \|^{\frac{3}{5}}_{L^2(0,T;L^6(\Omega))}
	  \| \nabla_{q_j} \phi \|_{L^\infty(0,T;L^\infty(\Omega \times D))}, 	
	\end{split}
	\end{equation*}
and therefore, by invoking \eqref{cc4} and \eqref{L61_11}$_2$, we deduce that $I_3$ converges to zero as $\ep \to 0$.
Similarly, by noting that \eqref{L61_10} implies that
$\int_D M \widehat\Psi \vec{q}_i \otimes \nabla_{q_i} \phi{\rm\,d}\vec{q} \in L^2(\Omega \times (0,T))^{3 \times 3}$,
$i=1,\dots,K$, we get by \eqref{conv_2} that
	\begin{equation*}
	I_4 \to \int_0^T \int_{\Omega \times D}
	[M (\nabla_x \vU ) \vec{q}_i ] \widehat\Psi \cdot \nabla_{q_j} \phi {\rm\,d}\vec{q}\dx\dt \quad \mbox{ as }\ep \to 0.
	\end{equation*}
Therefore the last term of \eqref{T61_3} converges to the last term in  \eqref{weak_lim_FP}.
Analogously, the third term of \eqref{T61_3} converges to the third term in  \eqref{weak_lim_FP}. In this way we obtain \eqref{weak_lim_FP} for smooth test functions
$\phi \in C([0,T]; C^\infty(\overline{\Omega \times D}))$. In order to extend the class of test functions for the limiting equation we use the density of the function space $C([0,T]; C^\infty(\overline{\Omega \times D}))$ in $L^2(0,T; H^s({\Omega \times D}))$, the embedding $H^s(\Omega \times D) \subset W^{1,\infty} (\Omega \times D)$,
\eqref{lem_46}, \eqref{cc4}, \eqref{ub_p_2}, \eqref{ub_p_2_2}, \eqref{L61_10}, and \eqref{conv_2}.

\subsection{Convergence of the convective part of the momentum equation}
In this section we concentrate on proving \eqref{co_1}. The proof is based on performing a Helmholtz decomposition of the momentum $\rho_\ep \vu_\ep= \vec{H} [\rho_\ep \vu_\ep] + \vec{H}^\perp [\rho_\ep \vu_\ep]$, the proof of the compactness of the solenoidal part in the decomposition, and the analysis of the acoustic equation governing the time-evolution of the gradient component in the decomposition.

Motivated by the discussion in Section 5.4.2 of \cite{FN}, we begin by decomposing the convective term as follows:
	\begin{equation}\label{AAc_0}
	\rho_\ep \vu_\ep \otimes \vu_\ep
	=
	\vec{H}[\rho_\ep \vu_\ep] \otimes \vu_\ep
	+\vec{H}^\perp[\rho_\ep \vu_\ep] \otimes \vec{H}[ \vu_\ep]
	+\vec{H}^\perp[\rho_\ep \vu_\ep] \otimes \vec{H}^\perp[ \vu_\ep] .
	\end{equation}
Let us emphasise that the solenoidal part of the momentum $\vec{H}[\rho_\ep \vu_\ep]$ does not exhibit oscillations in time and in particular it converges a.e. on the set $\Omega \times (0,T)$. In order to see this, we choose a test function in the momentum equation \eqref{T61_2} of the form $\vec{H}[\vec{w}]$, where
$\vec{w} \in C^\infty_c(\overline{\Omega}\times [0,T))^3$, $\vec{w} \cdot \vec{n} |_{\partial \Omega} = 0$. Thanks to the uniform estimates obtained in Section~\ref{ss:uni}, by noting that the singular term is irrelevant as $\divx \vec{H}[\vec{w}] = 0$, we deduce that the family  $\{ \int_\Omega  \vec{H}[\rho_\ep \vu_\ep ] \cdot \vec{w} \dx \}_{\ep>0} $ forms a bounded and equicontinuous sequence in $C([0,T])$. Therefore by the Arzel\`a--Ascoli theorem we have that
	\begin{equation*}
	\int_\Omega \vec{H}[\rho_\ep \vu_\ep] \cdot \vec{w} \dx
	\to
	 \int_\Omega \vec{H}[\bar\rho \vU] \cdot \vec{w} \dx \quad \mbox{ in } C([0,T])
	\end{equation*}
for any $\vec{w}$ as above. By a density argument and noting \eqref{conv_3_0}, we infer that
	\begin{equation}\label{AAc_1}
	\vec{H}[\rho_\ep \vu_\ep ] \to \vec{H}[\bar\rho \vU ] \quad \mbox{ in }C([0,T] ; L^{\frac{2\gamma}{\gamma+1}}_{weak}(\Omega)^3) .
	\end{equation}	
Next, we have that
	$$\bar\rho \vec{H}[\vu_\ep] \cdot \vu_\ep
	= \vec{H}[(\bar\rho - \rho_\ep)\vu_\ep] \cdot \vu_\ep
	+ \vec{H}[\rho_\ep \vu_\ep ] \cdot \vu_\ep.
	$$
The first term on the right-hand side of this equality converges to zero weakly in $L^1(\Omega \times (0,T))$ by \eqref{conv_1}, \eqref{conv_2} and compact  embedding. The second term, by \eqref{conv_2}, compact embedding, \eqref{AAc_1}, and since $\gamma > \frac{3}{2}$, converges weakly in $L^1(\Omega \times (0,T))$ to $\bar\rho|\vU|^2$. Hence, as $\vU = \vec{H}[\vU]$, we deduce that
	\begin{equation}\label{AAc_2}
	\vec{H}[\vu_\ep] \to \vU \quad \mbox{ strongly in }L^2(0,T; L^2(\Omega)^3).
	\end{equation}
This, together with \eqref{AAc_1},  implies that
	\begin{equation}\label{AAc_3}
	\vec{H}[\rho_\ep \vu_\ep] \otimes \vu_\ep \weak \bar\rho \vU \otimes \vU
	\quad \mbox{ weakly in } L^2(0,T;L^{\frac{6\gamma}{4\gamma +3}}(\Omega)^{3 \times 3}) .
	\end{equation}
Combining \eqref{AAc_2} and \eqref{conv_3} we get
	\begin{equation}\label{AAc_4}
	\vec{H}^\perp[\rho_\ep \vu_\ep] \otimes \vec{H} [\vu_\ep] \to \vec{0} \quad
	\mbox{ weakly in }  L^2(0,T;L^{\frac{6\gamma}{4\gamma +3}}(\Omega)^{3 \times 3}) .
	\end{equation}
Summarising \eqref{AAc_0}, \eqref{AAc_3}, and \eqref{AAc_4}, we see that in order to prove \eqref{co_1} all that is left to be shown is that
	\begin{equation}\label{AAc_5}
	\int_0^T \int_\Omega \vec{H}^\perp[\rho_\ep \vu_\ep] \otimes \vec{H}^\perp [\vu_\ep] : \nabla_x \vec{w} \dx\dt \to 0 \quad \mbox{ for any }\vec{w} \in C^\infty_c(\overline{\Omega} \times [0,T))^3, \ \divx \vec{w} = 0,\
	\vec{w} \cdot \vec{n} |_{\partial \Omega} = 0.
	\end{equation}
As our a-priori estimates do not provide any bound on the time-derivative of the gradient part of the velocity/momentum, in order to prove \eqref{AAc_5} we follow the strategy from \cite{FN}, which is based on the observation that possible oscillations in time mutually cancel in the acoustic waves described by means of $\vec{H}^\perp [\rho_\ep \vu_\ep]$ governed by the acoustic equation associated with our \eqref{NSFP_2} system.

\subsubsection{Acoustic equation}
The analysis of the following acoustic equation will allow us to control the temporal oscillations of the gradient part
of the momentum. Let us set, to this end,
	\begin{equation}
	\vec{V}_\ep := \re \ue, \quad \quad r_\ep := \frac{\re -\rs}{\ep},
	\end{equation}
and
	\begin{equation}
	\tens{L} := \tS (\ue ) - \re \ue \otimes \ue - \frac{1}{\ep^2} \left(p(\re) - p'(\rs) (\re - \rs) - p(\rs)\right)\tens{I} + \tens{\tau}(\psi_\ep).
	\end{equation}
We can then rewrite the continuity equation and the momentum equation with the extra stress tensor on its right-hand side in the form of Lighthill's acoustic analogy (cf. eqs. (5.137) and (5.140) in \cite{FN}):
	\begin{equation}\label{wave_1}
	\begin{split}
	\ep \partial_t r_\ep + \divx \vec{V}_\ep & = 0 \quad\quad \mbox{ in } (0,T) \times \Omega,\\
	\ep \partial_t  \vec{V}_\ep + p'(\rs) \nabla_x r_\ep & = \ep \divx \tens{L}  \quad \quad \mbox{ in } (0,T) \times \Omega,
	\end{split}
	\end{equation}
supplemented with the boundary condition 	
	\begin{equation}\label{wave_1bc}
	 \vec{V}_\ep \cdot \vec{n} |_{\partial \Omega} = 0 .
	 \end{equation}
The equations \eqref{wave_1}, together with boundary condition \eqref{wave_1bc}, are understood in a weak sense; specifically, the integral identity
	\begin{equation}\label{wave_2}
	\int_0^T \int_\Omega
	\left(
	\ep r_\ep \partial_t \varphi+  \vec{V}_\ep \cdot \nabla_x  \varphi
	\right) \dx\, \mathrm{d}t = 0
	\end{equation}
holds for any $\varphi \in C_c^\infty ( \overline\Omega  \times (0,T)   )$ and
	\begin{equation}\label{wave_3}
	\int_0^T \int_\Omega
	\left(
	\ep \vec{V}_\ep \cdot \partial_t \vec \varphi    + p'(\bar\rho) r_\ep \divx \vec\varphi
	\right) \dx\, \mathrm{d}t = \ep \int_0^T \int_\Omega  \tens{L}_\ep : \nabla_x \vec\varphi \dx\, \mathrm{d}t
	\end{equation}
holds for any $\vec\varphi \in C_c^\infty ( \overline\Omega \times  (0,T) )^3$, $\vec{\varphi} \cdot  \vec{n} |_{\partial\Omega} = 0$.	

By the uniform estimates \eqref{ue_1}, \eqref{ue_6} we obtain that
	\begin{equation}\label{ae_1}
	\| r_\ep \|_{L^\infty(0,T; (L^2 + L^q)(\Omega))} \leq c \quad \mbox{ for } q= \min\{\gamma,2\},
	\end{equation}
and by \eqref{ue_7}, \eqref{ue_9} we deduce that
	\begin{equation}\label{ae_2}
	\| \vec{V}_\ep \|_{L^\infty (0,T; (L^2 + L^s)(\Omega))} \leq c \quad \mbox{ with } s\in[1,2\gamma / (\gamma+1)] .
	\end{equation}
Moreover, by combining \eqref{ue_1}, \eqref{ue_5}, \eqref{ue_6}, \eqref{ue_8}, \eqref{L61_10}, \eqref{p428b} we deduce that
	$$ \tens{L} = \tens{L}_1 + \tens{L}_2 + \tens{L}_3 ,$$
where
	\begin{equation}\label{ae_3}
	\| \tens{L}_1 \|_{L^\infty(0,T; L^1(\Omega))} \leq c ,
	\quad
	 \| \tens{L}_2 \|_{L^2(0,T; L^2(\Omega))} \leq c, \quad \mbox{ and } \quad
	 \quad
	 \| \tens{L}_3\|_{L^2(0,T; L^{\frac{4}{3}}(\Omega))} \leq c.
	\end{equation}

Now we can follow the strategy developed in \cite{FN}. We reduce the problem to a finite number of modes, which are represented by eigenfunctions of the wave operator in \eqref{wave_2}, \eqref{wave_3}. It then transpires that the nonvanishing oscillatory part of the convective term can be represented as the gradient of a scalar function, and it can be therefore incorporated into the limiting pressure $\Pi$, whereby it is irrelevant in the incompressible limit of $\ep \rightarrow 0$. The details of the spectral analysis of the wave operator appearing in \eqref{wave_2}, \eqref{wave_3} are contained in the next section.

\subsubsection{Spectral analysis of the wave operator}
Next, as in Section 5.4.5 of \cite{FN}, we now focus on the following eigenvalue problem associated with the operator appearing on the left-hand side of \eqref{wave_2}, \eqref{wave_3}:
	\begin{equation}\label{wo_1}
	\divx \vec{\phi} = - \lambda \zeta, \quad p'(\bar\rho) \nabla_x \zeta = \lambda \vec{\phi} \ \mbox{ in } \Omega, \quad \vec{\phi} \cdot \vec{n}|_{\partial \Omega} = 0,
	\end{equation}
and reformulate it as the following homogeneous Neumann problem:
	\begin{equation}\label{wo_2}
	- \Delta_x \zeta = \Lambda \zeta \ \mbox{ in } \ \Omega,
	\quad \nabla_x \zeta \cdot \vec{n}|_{\partial \Omega}=0, \quad -\Lambda = \frac{\lambda^2}{p'(\bar\rho)},
	\end{equation}
for which we have a countable system of eigenvalues
$0=\Lambda_0 < \Lambda_1 \leq \Lambda_2\leq \cdots$, with the associated system of real eigenfunctions $\{ \zeta_n \}_{n=0}^\infty$ being a basis of the Hilbert space $L^2(\Omega)$. The complex eigenfunctions $\{ \vec{\phi}_{\pm n} \}_{n=0}^\infty$ are determined through  \eqref{wo_1} as
	\begin{equation}\label{wo_3}
	 \vec{\phi}_{\pm n}  := \pm \sqrt{\frac{p'(\bar\rho)}{\Lambda_n}} \nabla_x \zeta_n, \quad n=1,2,\dots .
	\end{equation}
Moreover, we can decompose the Hilbert space $L^2(\Omega)^3= L^2_{\mathrm{div}}(\Omega)^3 \oplus H_{\perp} (\Omega)$, where $H_\perp(\Omega)$ is the closure in the $L^2(\Omega)^3$ norm of the space spanned by $\{( {-i}/{\sqrt{p'(\bar\rho)}}) \vec{\phi}_n \}_{n=1}^\infty $, and where  $L^2_{\mathrm{div}}(\Omega)^3$ denotes the subspace of divergence-free  3-component vector functions contained in $L^2(\Omega)^3$, with vanishing normal trace on $\partial\Omega$; that is, the closure in the $L^2(\Omega)^3$-norm of the set of all $\vec{\varphi} \in C_c^\infty(\Omega)^3$ s.t. $\divx \vec{\varphi} = 0$ in $\Omega$. In order to reduce the problem to a finite number of modes, we introduce the corresponding orthogonal projection
	$$
	\vec{P}_N : L^2(\Omega)^3
	\to {\rm span}\,
	\bigg\{ \frac{-i}{\sqrt{p'(\bar\rho)}} \vec{\phi}_n \bigg\}_{n \leq N}, \quad N = 1,2,\dots
	$$
and we set
	$$\vec{H}^\perp_N [\vec\varphi] := \vec{P}_N \vec{H}^\perp[\vec\varphi]
	=  \vec{H}^\perp  \vec{P}_N[\vec\varphi]  $$
(since $\vec{P}_N$ commutes with $\vec{H}^\perp$).
For any $\vec\varphi \in L^2(\Omega)^3$ we consider the Fourier coefficients
	\begin{equation}\label{wo_3_0}
	a_n [\vec\varphi] := \frac{-i}{\sqrt{p'(\bar\rho)}}
	\int_\Omega \vec\varphi  \cdot \vec\phi_n \dx
	\end{equation}
and a scale of Hilbert spaces $H_{\alpha,\perp}(\Omega) \subset H_{\perp}(\Omega)$, $\alpha \in [0,1]$, with the norm $\| \cdot \|_{L^2_{\alpha,\perp}}$, defined by $$\| \vec\varphi \|^2_{L^2_{\alpha,\perp}} := \sum_{n=1}^\infty \Lambda^\alpha_n\, | a_n [\vec\varphi]|^2,$$
where $\{\Lambda_n\}_{n=1}^\infty \subset \mathbb{R}_{>0}$
is the family of nonzero eigenvalues associated with \eqref{wo_2}.
We note that
	\begin{equation}\label{wo_4}
	\begin{split}
	\| \vec{H}^\perp[\vec\varphi] -  \vec{H}^\perp_N[\vec\varphi] \|^2_{L^2_{\alpha_1,\perp}}
	& = \sum_{n=N+1}^\infty \Lambda_n^{\alpha_1} | a_n[\vec\varphi]|^2
	\leq \Lambda_N^{\alpha_1 - \alpha_2}  \sum_{n=N+1}^\infty \Lambda_n^{\alpha_2} | a_n[\vec\varphi]|^2
	\\ &
	=  \Lambda_N^{\alpha_1- \alpha_2}
	 \|\vec{H}^\perp [\vec\varphi] - \vec{H}^\perp_N [\vec\varphi]\|^2_{L^2_{\alpha_2,\perp}}
	  \quad \mbox{ for } \alpha_2 > \alpha_1\,\, \mbox{and any integer $N \geq 0$}.
	\end{split}
	\end{equation}
As $H_{0,\perp} = H_{\perp}$ and $H_{1,\perp} \subset L^6(\Omega)^3$, we deduce by an interpolation  argument that there exists an  $\bar\alpha \in (0,1)$ such that
	\begin{equation}\label{wo_5}
	H_{\alpha,\perp}(\Omega) \subset L^{s'}(\Omega)^3 \quad \mbox{ where }\, s = \frac{2\gamma}{\gamma +1}\; \mbox{ and }\; s'=\frac{2\gamma}{\gamma -1}\;
	\mbox{ whenever } \;\alpha \geq \bar\alpha.
	\end{equation}
Let us return to \eqref{AAc_5} and rewrite it as follows:
	\begin{equation*}
	\begin{split}
	\int_0^T \int_\Omega &
	 \vec{H}^\perp[\rho_\ep \vu_\ep] \otimes \vec{H}^\perp [\vu_\ep] : \nabla_x \vec{w} \dx\dt
	 \\ & =
	 \int_0^T \int_\Omega \vec{H}^\perp[\rho_\ep \vu_\ep] \otimes \vec{H}^\perp_N [\vu_\ep] : \nabla_x \vec{w} \dx\dt
	 + \int_0^T \int_\Omega \vec{H}^\perp[\rho_\ep \vu_\ep] \otimes ( \vec{H}^\perp [\vu_\ep]  - \vec{H}^\perp_N [\vu_\ep] ): \nabla_x \vec{w} \dx\dt .
	\end{split}
	\end{equation*}
 By \eqref{conv_3_0}, \eqref{wo_4} and \eqref{wo_5} we deduce that
 	$$
	\left|
	 \int_0^T \int_\Omega \vec{H}^\perp[\rho_\ep \vu_\ep] \otimes ( \vec{H}^\perp [\vu_\ep]  - \vec{H}^\perp_N [\vu_\ep] ): \nabla_x \vec{w} \dx\dt
	 \right| \leq c \Lambda_N^{-\frac{1}{2}(1-\bar\alpha)},
	$$
uniformly as $\ep \to 0$, for any fixed $\vec{w}
\in C^\infty_c(\overline{\Omega} \times [0,T))^3$, such that $\divx \vec{w} = 0$ and  $\vec{w} \cdot \vec{n} |_{\partial \Omega} = 0$.

We note that $\Lambda_N \to \infty$ as $N \to \infty$.
By a duality argument and \eqref{wo_5} we have also that
	$$
	\| \vec{H}^\perp [\vec\varphi]  - \vec{H}^\perp_N [\vec\varphi] \|^2_{[H^1(\Omega)]^*}
	 \leq c \Lambda_N^{\bar\alpha-1}
	 \| \vec\varphi \|^2_{L^{\frac{2\gamma}{\gamma +1}}(\Omega)}.
	$$
Indeed, this bound is an immediate consequence of the following:
\begin{align*}
\int_\Omega (\vec{H}^\perp [\vec\varphi]  - \vec{H}^\perp_N [\vec\varphi])\cdot  \vec{\psi}  \dx 
&
= \int_\Omega \vec\varphi \cdot  (\vec{H}^\perp [\vec\psi]  - \vec{H}^\perp_N [\vec\psi])\dx
\\
& 
\leq \|\vec\varphi\|_{L^{\frac{2\gamma}{\gamma +1}}(\Omega)} \|\vec{H}^\perp [\vec\psi]  - \vec{H}^\perp_N [\vec\psi]\|_{L^{\frac{2\gamma}{\gamma-1}}(\Omega)}\\
& \leq \|\vec\varphi\|_{L^{\frac{2\gamma}{\gamma +1}}(\Omega)} \| \vec{H}^\perp[\vec\psi] -  \vec{H}^\perp_N[\vec\psi] \|_{L^2_{\bar\alpha,\perp}}\\
& \leq \|\vec\varphi\|_{L^{\frac{2\gamma}{\gamma +1}}(\Omega)} \Lambda_N^{\frac{1}{2}(\bar\alpha- 1)}
	 \|\vec{H}^\perp [\vec\psi] - \vec{H}^\perp_N [\vec\psi]\|_{L^2_{1,\perp}}\\
& \leq c \|\vec\varphi\|_{L^{\frac{2\gamma}{\gamma +1}}(\Omega)} \Lambda_N^{\frac{1}{2}(\bar\alpha- 1)}
	 (\|\vec{H}^\perp [\vec\psi]\|_{L^2_{1,\perp}} + \|\vec{H}^\perp_N [\vec\psi]\|_{L^2_{1,\perp}})\\
&\leq c \|\vec\varphi\|_{L^{\frac{2\gamma}{\gamma +1}}(\Omega)} \Lambda_N^{\frac{1}{2}(\bar\alpha- 1)}
	 (\|\vec\psi\|_{L^2_{1,\perp}} + \|\vec\psi\|_{L^2_{1,\perp}})\\
& \leq c \|\vec\varphi\|_{L^{\frac{2\gamma}{\gamma +1}}(\Omega)} \Lambda_N^{\frac{1}{2}(\bar\alpha- 1)}
	 \|\vec\psi\|_{H^1(\Omega)},
\end{align*}
where the penultimate inequality in this chain of inequalities follows by noting that 
$a_n[\vec{H}^\perp [\vec\psi]]=a_n[\vec\psi]$ for all $n \geq 1$, that $a_n[\vec{H}^\perp_N [\vec\psi]] = a_n[\vec\psi]$ for $1 \leq n \leq N$ and $a_n[\vec{H}^\perp_N [\vec\psi]] = 0$ for $n > N$, whereby $\|\vec{H}^\perp [\vec\psi]\|_{L^2_{1,\perp}} = \|\vec\psi\|_{L^2_{1,\perp}}$ and
$\|\vec{H}^\perp_N [\vec\psi]\|_{L^2_{1,\perp}} \leq  \|\vec\psi\|_{L^2_{1,\perp}}$.

With these arguments in hand, the proof of \eqref{AAc_5} reduces to showing that
	$$
	\int_0^T \int_\Omega
	 \vec{H}^\perp_N[\rho_\ep \vu_\ep] \otimes \vec{H}^\perp_N [\vu_\ep] : \nabla_x \vec{w} \dx\dt \to 0\mbox{ for any }\vec{w} \in C^\infty_c(\overline{\Omega} \times [0,T))^3, \ \divx \vec{w} = 0,\
	\vec{w} \cdot \vec{n} |_{\partial \Omega} = 0
	$$
as $\ep \to 0$, or, rather, thanks to \eqref{conv_1}, it suffices to show that
	\begin{equation}\label{wo_6}
	\int_0^T \int_\Omega
	 \vec{H}^\perp_N[\rho_\ep \vu_\ep] \otimes \vec{H}^\perp_N [\rho_\ep\vu_\ep] : \nabla_x \vec{w} \dx\dt \to 0 \mbox{ for any }\vec{w} \in C^\infty_c(\overline{\Omega} \times [0,T))^3, \ \divx \vec{w} = 0,\
	\vec{w} \cdot \vec{n} |_{\partial \Omega} = 0
	\end{equation}
as $\ep \to 0$; this is done in the next subsection.

\subsection{Weak limit in the convective term}
In order to prove \eqref{wo_6} let us choose a test function for \eqref{wave_2} of the following `separated' form
	$$\varphi(\vec{x},t) = \kappa(t)\, \zeta_n(\vec{x}) \quad \mbox{ with } \kappa \in C^\infty_c (0,T),$$
where $\zeta_n$ and $\Lambda_n$ solve the eigenvalue problem \eqref{wo_2}, and for \eqref{wave_3}
as a test function we take
	$$\vec\varphi (\vec{x},t) = \kappa(t)\, \frac{1}{\sqrt{\Lambda_n}} \nabla_x \zeta_n \quad \mbox{ with } \kappa \in C^\infty_c(0,T) .$$
Consequently, we obtain the following system of ordinary differential equations:
	\begin{equation}\label{wo_7}
	\begin{split}
	\ep \partial_t( b_n[r_\ep]) - \sqrt{\Lambda_n}\, a_n[\vec{V}_\ep] & = 0 \quad \mbox{ in }  (0,T),
	\\
	\ep \partial_t( a_n[\vec{V}_\ep]) + p'(\bar\rho) \sqrt{\Lambda_n}\, b_n [r_\ep] & =  \ep L_{\ep,n} \quad \mbox{ in }  (0,T),
	\end{split}
	\end{equation}
for $n= 1,2, \dots $, where $a_n[\vec{V}_\ep]$ are the Fourier coefficients of $\vec{V}_\ep$, defined by \eqref{wo_3_0},
	$
	b_n[r_\ep] := \int_{\Omega} r_\ep \zeta_n \dx ,
	$
and
	\begin{equation}\label{wo_8}
	\| L_{\ep,n} \|_{L^1(0,T)} \leq c\quad \mbox{ for any fixed }n=1,2,3,\dots .
	\end{equation}
Then, in terms of the Helmholtz projection, \eqref{wo_7} reads as follows:
	\begin{equation}\label{wo_9}
	\begin{split}
	\ep \partial_t ([r_\ep]_N) + \divx (\vec{H}^\perp_N[\rho_\ep \vu_\ep])  & = 0 \quad \mbox{ in }  \Omega \times (0,T),
	\\
	\ep \partial_t (\vec{H}^\perp [\rho_\ep \vu_\ep]) + p'(\bar\rho) \nabla_x [r_\ep]
	& = \ep  L_{\ep,N} \quad \mbox{ in }  (0,T) \times \Omega,
	\end{split}
	\end{equation}
where we set $[r_\ep]_N := \sum_{n=1}^N b_n [r_\ep] \zeta_n$,
and thanks to \eqref{wo_8} we get that
$\| L_{\ep,N}\|_{L^1(\Omega \times (0,T))} \leq  c$ for any fixed $N\geq 1$. Notice that \eqref{wo_9} is satisfied in a strong sense, since $[r_\ep]_N$ and $\vec{H}^\perp [\rho_\ep \vu_\ep]$ are regular enough. By introducing the potential $\Phi_{\ep,N}$ via
	$$\nabla_x \Phi_{\ep,N} = \vec{H}^\perp_N[\rho_\ep \vu_\ep], \quad \int_\Omega\Phi_{\ep,N} \dx = 0,$$
we can reformulate \eqref{wo_6} as
	\begin{equation}\label{wo_10}
	\int_0^T \int_\Omega \vec{H}^\perp_N [\rho_\ep \vu_\ep] \otimes \vec{H}^\perp  [\rho_\ep \vu_\ep]  : \nabla_x \vec{w} \dx\dt = - \int_0^T \int_\Omega \Delta_x \Phi_{\ep,N} \nabla_x \Phi_{\ep,N} \cdot \vec{w} \dx\dt,
	\end{equation}
where $\vec{w} \in C^\infty_c(\overline{\Omega} \times (0,T))^3$, $\divx \vec{w} = 0$,
$\vec{w} \cdot \vec{n} |_{\partial \Omega}= 0.$
Then by \eqref{wo_9} we obtain that
	\begin{equation}\label{wo_11}
	\int_0^T \int_\Omega \Delta_x \Phi_{\ep,N} \nabla_x \Phi_{\ep,N} \cdot \vec{w} \dx\dt
	= \ep \int_0^T \int_\Omega [r_\ep]_N \nabla_x \Phi_{\ep,N} \cdot \partial_t \vec{w} \dx\dt
	+ \ep \int_{0}^T \int_\Omega [r_\ep]_N L_{\ep,N} \cdot \vec{w} \dx\, \mathrm{d}t .
	\end{equation} 	
Therefore, thanks to \eqref{wo_8}, the  right-hand side of \eqref{wo_11} converges to zero as $\ep \to 0$ for any fixed $\vec{w}$ as above.  This step completes the proof of \eqref{wo_6} and consequently we may replace $\overline{\bar\rho \vU \otimes \vU}$ by  ${\bar\rho \vU \otimes \vU}$ in \eqref{weak_lim_mom_0}, and deduce that \eqref{weak_lim_mom} is satisfied.

The formulation \eqref{weak_lim_mom} with the incompressibility constraint \eqref{divU} is supplemented by the boundary conditions
	\begin{equation}
	\vU \cdot \vec{n}|_{\partial \Omega} = 0, \quad [\tens{D}\vU]\vec{n} \times \vec{n}|_{\partial \Omega} = \vec{0}.
	\end{equation}
Moreover, using \eqref{AAc_1} we can show that $\vU \in C([0,T]; L^2_{weak}(\Omega))$, and we thus deduce also
that
	$$
	\int_\Omega \vU(0,\cdot) \cdot \vec{w} \dx = \int_\Omega \vU_0 \cdot \vec{w} \dx  \quad
	\mbox{ for all } \vec{w} \in C^\infty_c(\overline{\Omega}), \ \divx \vec{w} = 0; \
\vec{w} \cdot \vec{n} |_{\partial \Omega}= 0;
	 $$
in other words, $\vU(0,\cdot) = \vec{H}[\vU_0]$.

\bigskip

The formulation \eqref{weak_lim_mom} is satisfied for all functions $\vec{w}$ that are smooth, divergence-free, and $\vec{w} \cdot{\vec{n}}|_{\partial \Omega} = 0$. However this family of test functions can be extended to the one mentioned in \eqref{weak_lim_mom} by a density argument and the estimates obtained in Section~\ref{ss:uni}.

\section{Concluding remarks}\label{sec3}
We studied the behaviour of global-in-time weak solutions to a class of bead-spring chain models with finitely extensible nonlinear elastic (FENE) spring potentials for dilute polymeric fluids, and we proved that as the Mach number tends to zero the system is driven to its incompressible counterpart. Our analysis was performed under the assumption that the velocity of the solvent in which the polymer molecules are immersed satisfies a complete slip boundary condition. The corresponding passage to the limit in the case of a nonslip boundary condition for the velocity field is, as in the case in the compressible Navier--Stokes system (i.e. in the absence of coupling to the Fokker--Planck equation), more technical (cf. Ch.7 of \cite{FN}, particularly the first paragraph of subsection 7.1.2) and will be considered elsewhere. The existence of global weak solutions to the coupled incompressible Navier--Stokes--Fokker--Planck system was proved in \cite{BS10} on arbitrary bounded open Lipschitz domains $\Omega \subset
\mathbb{R}^d$, $d \in \{2,3\}$, while the existence proof in \cite{BS2014} for the corresponding compressible Navier--Stokes--Fokker--Planck system was restricted to bounded domains $\Omega \subset \mathbb{R}^d$, $d \in \{2,3\}$, with $\partial\Omega \subset C^{2,\alpha}$, $\alpha \in (0,1)$. Using the ideas in \cite{FNP} the $C^{2,\alpha}$ regularity of $\partial\Omega$ assumed in \cite{BS2014} can be relaxed to $\Omega$ being a bounded open Lipschitz domain. For simplicity and for the sake of consistency with the assumptions on $\Omega$ in \cite{FN} we have however restricted ourselves here to domains of class $C^{2,\alpha}$.

\section*{Acknowledgments}
 AWK was partially supported by a Newton Fellowship of the Royal Society and by the grant Iuventus Plus 0871/IP3/2016/74 of Ministry of Sciences and Higher Education RP.

\bigskip
Aneta Wr\'oblewska-Kami\'nska\\
Institute of Mathematics\\
 Polish Academy of Sciences\\
  ul. \'Sniadeckich 8
  00-656 Warsaw, Poland
 \\[1.5ex]
Endre S\"uli\\
Mathematical Institute\\
University of Oxford\\
Woodstock Road\\
Oxford OX2 6GG, UK
\bigskip

\end{document}